\documentclass[]{amsart}

\usepackage{amsbsy, amsmath, amssymb, amstext, amsxtra, commath, eucal, float, graphicx, hyperref, latexsym, mathrsfs, mathtools, multicol, tikz, tikz-cd}

\usepackage{lineno}

\newtheorem{theorem}{Theorem}[section]
\newtheorem{lemma}[theorem]{Lemma}
\newtheorem{proposition}[theorem]{Proposition}
\newtheorem{corollary}[theorem]{Corollary}
\theoremstyle{definition}
\newtheorem{definition}[theorem]{Definition}

\newcommand\restr[2]{\ensuremath{#1\!\!\upharpoonright_{#2}}\,}

\def\A{\mathop{\mathsf{A}\hspace{0mm}}\nolimits}
\def\C{\mathop{\mathsf{C}\hspace{0mm}}\nolimits}
\def\cf{\mathop{\rm cf}\nolimits}
\def\CN{\mathop{\mathsf{CN}}\nolimits}
\def\int{\mathop{\rm int}\nolimits}
\def\log{\mathop{\rm log}\nolimits}
\def\ON{\mathop{\mathsf{ON}}\nolimits}
\def\P{\mathop{\mathsf{P}\hspace{0mm}}\nolimits}
\def\UR{\mathop{\mathsf{UR}}\nolimits}
\def\UC{\mathop{\mathsf{UC}}\nolimits}
\def\WP{\mathop{\mathsf{WP}}\nolimits}

\numberwithin{equation}{section}

\let\mathcal\mathscr

\begin{document}

%\linenumbers

\title{On calibers for $C_p(X)$}

\author[A. Ríos-Herrejón]{Alejandro R\'ios Herrej\'on}
\address {A. Ríos Herrej\'{o}n\\
Departamento de Matem\'aticas, Facultad de Ciencias, Universidad Nacional Aut\'onoma de M\'exico, Circuito ext. s/n, Ciudad Universitaria, C.P. 04510,  M\'exico, CDMX}
\email{chanchito@ciencias.unam.mx}
%\thanks{El primer autor recibió el apoyo económico del CONACYT (núm. 814282).}

\author[\'A. Tamariz-Mascar\'ua]{\'Angel Tamariz-Mascar\'ua}
\address{Á. Tamariz Mascar\'ua, Departamento de Matem\'aticas, Facultad de Ciencias, Universidad Nacional Aut\'onoma de M\'exico, Circuito ext. s/n, Ciudad Universitaria, C.P. 04510,  M\'exico, CDMX}
\email{atamariz@unam.mx}
\urladdr{}

\thanks{The first author was supported by CONAHCYT grant no. 814282.}



\begin{abstract} We present new results regarding calibers in the function spaces $C_p(X)$. Our main theorem is that $C_p(X)$ is strongly \v{S}anin whenever $X$ is a submetrizable space; this improves an earlier result due to Tkachuk: $C_p(X)$ is \v{S}anin whenever $X$ is a submetrizable space. Moreover, we give sufficient conditions to characterize the calibers of $C_p(X)$ when $X$ is a topological sum, and we calculate the calibers of $C_p(X)$ when $X = \prod_{\xi < \lambda}X_\xi$ is a product of non-trivial Tychonoff spaces with $i$-weight $\leq \lambda$. Furthermore, we calculate the calibers of $C_p(X)$ when $X$ is an interval of ordinals and when $X$ is the one-point $\lambda$-Lindelöf extension of a discrete space of cardinality $\geq \lambda$. This allows to give examples of compact Hausdorff spaces $Z$ such that $iw(Z)=\kappa^{+}$ and $C_p(Z^{\kappa})$ does not have caliber $iw(Z)$; and examples of spaces $\{Z_\alpha : \alpha<\cf(\kappa)\}$ such that $\kappa$ is a caliber for $C_p(Z_\alpha)$ whenever $\alpha<\cf(\kappa)$ but it is not a caliber for $C_p(\bigoplus_{\alpha<\cf(\kappa)} Z_\alpha)$.
\end{abstract}

\maketitle

%%%%%%%%%%%%%%%%%%%%%%%%%%%%%%%%%%%%%%%%%%%%%%%%%%%%%%%%%%%%%%
\section{Introduction}
%%%%%%%%%%%%%%%%%%%%%%%%%%%%%%%%%%%%%%%%%%%%%%%%%%%%%%%%%%%%%%

Since \v{S}anin presented the notion of {\it caliber} for a topological space in \cite{sanin1948}, this concept and its variants have been deeply studied. From classical texts like \cite{argtsa1982}, \cite{shelah1977} and \cite{tall1977}, to recent articles such as \cite{rios2022}, \cite{rios2023} and \cite{riotam2023}, chain conditions have proven to be a fruitful object of study.

In this paper we examine the chain conditions of the spaces of continuous real-valued functions defined on a Tychonoff topological space $X$, denoted by $C(X)$, with the pointwise convergence topology (in symbols: $C_p(X)$); that is, $C(X)$ considered as a subspace of $\mathbb{R}^{X}$. Our goal is to determine the calibers for $C_p(X)$ by means of the topological properties of $X$.

Regarding this problem, there is information in mathematical literature that is distributed in several texts (see, for example, \cite{kalspi1992} and \cite{tkachuk2015}--\cite{tkachuk1986}). In the following pages we make a synthetic compilation of some known theorems and prove several new results.

In Section~\ref{secc_preliminaries} we gather, without proof, some results we use throughout the text. In Section~\ref{secc_Cp(X)} we provide basic results about the calibers of $C_p(X)$. In Sections~\ref{secc_submetrizable}, \ref{secc_ordinales} and \ref{secc_Lindelof} we calculate the calibers of $C_p(X)$ where $X$ is, respectively, a submetrizable space, an interval of ordinals and the one-point $\lambda$-Lindelöf extension of a discrete space of cardinality $\geq \lambda$; in particular, we prove that $C_p(X)$ is strongly \v{S}anin whenever $X$ is submetrizable. Finally, in Sections~\ref{secc_productos} and~\ref{secc_sumas} we study the calibers of $C_p(X)$ when $X$ is, correspondingly, a topological product and a topological sum.

%%%%%%%%%%%%%%%%%%%%%%%%%%%%%%%%%%%%%%%%%%%%%%%%%%%%%%%%%%%%%%
\section{Preliminaries}\label{secc_preliminaries}
%%%%%%%%%%%%%%%%%%%%%%%%%%%%%%%%%%%%%%%%%%%%%%%%%%%%%%%%%%%%%%

Our reference text for basic topology material is \cite{engelking1989}. Any topological notion that is not explicitly defined in the following paragraphs should be understood as indicated there.

Throughout this article $\omega$ is the first infinite ordinal (hence, the first infinite cardinal) and $\mathbb{N} := \omega \setminus \{0\}$. Furthermore, a set $X$ is {\it countable} if there is an injective function from $X$ into $\omega$.

In this work we denote by $\ON$ and $\CN$ the proper classes formed by the ordinal numbers and the infinite cardinal numbers, respectively. Furthermore, we use the symbols $\UR$ and $\UC$ to refer to the subcollections of $\CN$ formed by the uncountable regular cardinals and the cardinals with uncountable cofinality, correspondingly.

If $X$ is a topological space, we denote by $\tau_{X}$ the topology of $X$ and by $\tau_{X}^{+}$ the set $\tau_{X}\setminus \{ \emptyset\}$. For a point $x\in X$, $\tau_{X}(x)$ denotes the subset of $\tau_{X}$ formed by the open sets that have $x$ as an element.

If $\mathcal{U}$ is a pairwise disjoint subset of $\tau_{X}^{+}$, we say that $\mathcal{U}$ is a \textit{cellular family} in $X$; and $X$ satisfies the \textit{countable chain condition} (abbreviated \textit{ccc}) if all cellular families in $X$ are countable.
The terms {\it centered family} and {\it family with the finite intersection property} in $X$ are used, indistinctly, to designate a collection of subsets such that the intersection of the elements of any of their finite subcollections is non-empty.
On the other hand, we say that a subset $\mathcal{U}$ of $\tau_X^+$ is {\it linked} if $U\cap V\neq \emptyset$ holds for every $U,V\in \mathcal{U}$.

If $\kappa$ is an infinite cardinal we say that $X$ {\it has precaliber $\kappa$} (resp., {\it weak precaliber $\kappa$}) if for any set $I$ of cardinality $\kappa$ and for any family $\{U_\alpha : \alpha\in I \} \subseteq \tau_{X}^{+}$, there exists $J\in [I]^{\kappa}$ such that $\{U_\alpha : \alpha\in J\}$ is a centered (resp., linked) family.
Moreover, we say that $\kappa$ is a {\it caliber} for $X$ if for every set $I$ of cardinality
$\kappa$ and for every family $\{U_\alpha : \alpha\in I \} \subseteq \tau_{X}^{+}$, there exists $J\in [I]^{\kappa}$ with $\bigcap_{\alpha\in J} U_\alpha \neq\emptyset$.

We work throughout this article with the following collections of cardinal numbers: \begin{align*} \WP(X)&:=\{\kappa\in\CN : \kappa \ \text{is a weak precaliber for} \ X\}, \\
\P(X)&:=\{\kappa\in\CN : \kappa \ \text{is a precaliber for} \ X\}, \ \text{and} \\ 
\C(X)&:=\{\kappa\in\CN : \kappa \ \text{is a caliber for} \ X\}.
\end{align*}

A topological space $X$ is a {\it \v{S}anin space} if $\UR\subseteq \C(X)$. Similarly, we say that $X$ is a {\it strongly \v{S}anin space} if $\UC\subseteq \C(X)$.

Additionally, we work with the following cardinal functions (see \cite{hodel1984}) which are called, respectively, the {\it weight}, the {\it $\pi$-weight}, the {\it net weight}, the {\it density}, the {\it extent} and the {\it cellularity} of $X$: \begin{align*} w(X) &:= \min\left\{|\mathcal{B}| : \mathcal{B} \ \text{is a base for} \ X \right\}+\omega, \\
\pi w(X) &:= \min\left\{|\mathcal{V}| : \mathcal{V} \ \text{is a $\pi$-base for} \ X \right\}+\omega, \\
nw(X) &:= \min\left\{|\mathcal{N}| : \mathcal{N} \ \text{is a net for} \ X \right\}+\omega, \\
d(X) &:= \min\left\{|D| : D \ \text{is a dense subset of} \ X \right\} +\omega, \\
e(X) &:= \sup\left\{|D| : D \ \text{is a closed and discrete subspace of} \ X\right\}+\omega, \ \text{and} \\
c(X) &:= \sup\left\{|\mathcal{U}| : \mathcal{U} \ \text{is a cellular family in}\ X \right\} +\omega.
\end{align*}

For a metrizable space $X$ with metric $d$, $x\in X$ and $\varepsilon>0$, the set $B_\varepsilon(x)$ is the {\it open ball with center $x$ and radius $\varepsilon$}, i.e., $B_\varepsilon(x) := \{y\in X : d(x,y)<\varepsilon\}$.

Below we present (without proof) some auxiliary well-known results which are valid for every topological space.

\begin{proposition}\label{prop_caliber_density} If $\kappa$ is a cardinal number, and $d(X)<\cf(\kappa)$, then $\kappa$ is a caliber for $X$. That is, $\{\kappa\in\CN : \cf(\kappa)>d(X)\}\subseteq \C(X)$.

\end{proposition}

\begin{proposition}\label{prop_caliber_cofinality} If $\kappa$ is a caliber (resp., precaliber; weak precaliber) for $X$, then $\cf(\kappa)$ is a caliber (resp., precaliber; weak precaliber) for $X$.

\end{proposition}

\begin{proposition}\label{prop_caliber_dense_subset} If $D$ is a dense subspace of $X$ and $\kappa$ is an infinite cardinal, then the following statements hold.
\begin{enumerate}
\item $c(D)=c(X)$.
\item $\kappa$ is a precaliber (resp., weak precaliber) for $D$ if and only if $\kappa$ is a precaliber (resp., weak precaliber) for $X$.
\item If $\kappa$ is a caliber for $D$, then $\kappa$ is a caliber for $X$.
\end{enumerate}

\end{proposition}

\begin{proposition}\label{cor_cellularity_cofinality} If $X$ is infinite and $T_2$, then $\WP(X)\subseteq \UC$.

\end{proposition}

\begin{proposition}\label{prop_caliber_pi_weight} If $\kappa$ is an infinite cardinal with $\pi w(X)<\kappa$ and $X$ has caliber $\cf(\kappa)$, then $\kappa$ is a caliber for $X$. In symbols, $\{\kappa\in\CN : \kappa>\pi w(X)\ \wedge \ \cf(\kappa)\in \C(X)\} \subseteq \C( X)$. Hence, if $X$ is a \v{S}anin space and $\pi w(X)<\aleph_{\omega_1}$, then $X$ is strongly \v{S}anin.

\end{proposition}

The next theorem compiles some of the results proved in \cite{rios2022}, \cite{sanin1948} and \cite{shelah1977}.

\begin{theorem}\label{thm_RSS} Let $\kappa$ be an infinite cardinal, $\lambda$ a cardinal number, $\{X_\alpha : \alpha<\lambda\}$ a family of topological spaces, $X$ a topological space and $\A\in \{\C,\P,\WP\}$.

\begin{enumerate}
\item $\UR\cap \bigcap_{\alpha<\lambda} \A(X_\alpha) \subseteq \A(\prod_{\alpha<\lambda} X_\alpha) \subseteq \bigcap_{\alpha<\lambda} \A(X_\alpha)$.

\item If $\lambda<\cf(\kappa)$ and $\kappa \in \bigcap_{\alpha<\lambda} \A(X_\alpha)$, then $\kappa\in \A(\prod_{\alpha<\lambda} X_\alpha)$.

\item $\UC\cap \A(X^\lambda) = \UC\cap \A(X)$.

\end{enumerate}

\end{theorem}

Finally, two combinatorial results that we will need below, whose proofs can be consulted in \cite[2.4, p.~9]{hodel1984} and \cite[Theorem~2.10, p.~5]{rios2022}, are:

\begin{theorem}\label{thm_regular_delta_system} If $\kappa$ is a cardinal such that $\kappa=\cf(\kappa)>\omega$ and $\{A_\alpha : \alpha<\kappa\}$ is a family of finite sets, then there are sets $J\in [\kappa]^{\kappa}$ and $A$ such that $A_\alpha \cap A_\beta = A$, whenever $\alpha, \beta \in J$ are distinct.

\end{theorem}

\begin{theorem}\label{thm_singular_delta_system} Let $\kappa$ be a cardinal and $\{\kappa_\alpha : \alpha<\cf(\kappa)\}$ a family of regular uncountable cardinals with the following properties: 
\begin{enumerate}
\item $\kappa>\cf(\kappa)>\omega$,
\item $\sup\{\kappa_\alpha : \alpha<\cf(\kappa)\}=\kappa$, and
\item $\max\{\cf(\kappa), \sup\{\kappa_\beta : \beta<\alpha\}\}<\kappa_\alpha$ for every $\alpha<\cf(\kappa)$.
\end{enumerate} If for each $\alpha<\cf(\kappa)$  the set $I_\alpha$ has cardinality $\kappa_\alpha$ and $\{A(\alpha, \gamma) : \alpha<\cf(\kappa)\ \wedge\ \gamma\in I_\alpha\}$ is a family of finite sets, then there are $m<\omega$, $I\in \left[\cf(\kappa)\right]^{\cf(\kappa)}$, $\{J_\alpha : \alpha\in I\}$, $\{A_\alpha : \alpha\in I\}$ and $A$ that satisfy the following conditions for different $\alpha,\beta\in I$: 
\begin{enumerate}
\item $J_\alpha \in [I_\alpha]^{\kappa_\alpha}$,
\item $A_\alpha \cap A_\beta = A$,
\item if $\gamma,\delta \in J_\alpha$ are different, then $A(\alpha, \gamma) \cap A(\alpha, \delta) = A_\alpha$, 
\item whenever $\gamma\in J_\alpha$ and $\delta\in J_\beta$, $A(\alpha, \gamma) \cap A(\beta, \delta) = A$, and
\item $|A_\alpha|=m$.
\end{enumerate}

\end{theorem}

%%%%%%%%%%%%%%%%%%%%%%%%%%%%%%%%%%%%%%%%%%%%%%%%%%%%%%%%%%%%%%
\section{Chain conditions in $C_p$-spaces}\label{secc_Cp(X)}
%%%%%%%%%%%%%%%%%%%%%%%%%%%%%%%%%%%%%%%%%%%%%%%%%%%%%%%%%%%%%%

From now on, $X$ will be a Tychonoff space and $C(X)$ will be the collection of continuous functions from $X$ to $\mathbb{R}$. For every $f\in C(X)$, $A\in [X]^{<\omega}$ and $\varepsilon>0$ set $$V(f,A,\varepsilon) := \left\{g\in C(X) : \forall x\in A \left(|f(x)-g(x)|<\varepsilon\right)\right\}.$$

We use the symbol $C_p(X)$ to refer to the set $C(X)$ equipped with the topology of pointwise convergence; in other words, the topology of $C_p(X)$ is generated by the collection $$\left\{V(f,A,\varepsilon) : f\in C(X) \ \wedge\ A\in [X]^{<\omega}\setminus \{\emptyset\}\ \wedge\ \varepsilon>0\right\}$$ as a base. It is easy to check that $C_p(X)$ is a dense subspace of $\mathbb{R}^{X}$ (see \cite[Proposition~3.3.6, p.~9]{arh1992}).

The study of some chain conditions in the spaces $\mathbb{R}^{X}$ and $C_p(X)$ is trivialized by the next result:

\begin{proposition}\label{prop_Cp_precalibers} The following statements are true. \begin{enumerate} \item $\C\left(\mathbb{R}^{X}\right)=\P\left(\mathbb{R}^{X}\right)=\WP\left(\mathbb{R}^{X}\right)=\UC$.
\item $\P\left(C_p(X)\right)=\WP\left(C_p(X)\right)=\UC$.

\end{enumerate}

\end{proposition}

\begin{proof} Since $\mathbb{R}$ is an infinite and separable Hausdorff space, Proposition~\ref{prop_caliber_density} and Proposition~\ref{cor_cellularity_cofinality} ensure that $\C\left(\mathbb{R} \right)=\P\left(\mathbb{R}\right)=\WP\left(\mathbb{R}\right)=\UC$. Thus, Theorem~\ref{thm_RSS} implies that $\C\left(\mathbb{R}^{X}\right)=\P\left(\mathbb{R}^{X}\right) =\WP\left(\mathbb{R}^{X}\right)=\UC$. Finally, the relations $\P\left(C_p(X)\right)=\WP\left(C_p(X)\right)=\UC$ follow from Proposition~\ref{prop_caliber_dense_subset}.
\end{proof}

Now, since $\C(C_p(X))\subseteq \P\left(C_p(X)\right)$ it follows that \begin{align} \C\left(C_p(X)\right )\subseteq \UC. \label{inclusion}
\end{align} We are interested, among other things, in finding out under what conditions the opposite inclusion is verified: \begin{align} \UC\subseteq \C(C_p(X));\label{inclusion_2}
\end{align} in other words, we want to know when $C_p(X)$ is strongly \v{S}anin. For example, the following result is a consequence of Proposition~\ref{prop_Cp_precalibers}:

\begin{proposition}\label{prop_Cp_calibers_discrete} If $X$ is a discrete space, then $C_p(X)$ is strongly \v{S}anin.

\end{proposition}

Under certain restrictions on some cardinal functions of $X$, we can guarantee that~(\ref{inclusion_2}) is satisfied. For our first results in this direction it is convenient to establish some auxiliary facts.

We say that a family of basic open sets in $C_p(X)$ is {\it elementary} if it is of the form $\{V(f_\alpha,A_\alpha,\varepsilon) : \alpha<\kappa\}$, where $\varepsilon>0$ and all $A_\alpha$ have the same finite cardinality. This definition becomes relevant due to the following result.

\begin{lemma}\label{lemma_elementary_family} The following statements are equivalent for a topological space $X$ and a cardinal number $\kappa$ with $\cf(\kappa)>\omega$.

\begin{enumerate}
\item $C_p(X)$ has caliber $\kappa$.
\item For every elementary family $\{V(f_\alpha,A_\alpha,\varepsilon) : \alpha<\kappa\}$ there are $J\in [\kappa]^{\kappa}$ and $f\in C_p(X)$ with $f\in \bigcap\{V(f_\alpha,A_\alpha,\varepsilon) : \alpha\in J\}$.
\end{enumerate}

\end{lemma}

\begin{proof} The non-trivial implication is (2) implies (1). Let $\{V(f_\alpha,A_\alpha,\varepsilon_\alpha) : \alpha<\kappa\}$ be a family of basic open sets in $C_p(X)$. Since $\mathbb{R}$ has caliber $\kappa$, by Proposition~\ref{prop_Cp_precalibers}(1) there are $J\in [\kappa]^{\kappa}$ and $\varepsilon \in \bigcap\{(0,\varepsilon_\alpha) : \alpha\in J\}$. Furthermore, since the function $J \to \omega$ determined by $\alpha \mapsto \abs{A_\alpha}$ has a cofinal fiber, there exists $M \in [J]^{\kappa}$ with $\abs {A_\alpha} = \abs{A_\beta}$ whenever $\alpha,\beta \in M$. Thus, since $\{V(f_\alpha,A_\alpha,\varepsilon) : \alpha\in M\}$ is elementary, there are $L\in [M]^{\kappa}$ and $f\in C_p(X)$ with $f \in \bigcap\{V(f_\alpha,A_\alpha,\varepsilon) : \alpha\in L\}$. Finally, $L\in [\kappa]^{\kappa}$ and $f \in \bigcap\{V(f_\alpha,A_\alpha,\varepsilon_\alpha) : \alpha\in L\}$.
\end{proof}

If $f:X\to Y$ is a continuous function, the function $f^{\#}: C(Y) \to C(X)$ defined by $f^{\#}(g ) := g\circ f$ is called the {\it dual function} of $f$. In particular, when $f:X\to Y$ is a condensation between Tychonoff spaces, it is known that $f^{\#}: C_p(Y) \to C_p(X)$ is a topological embedding with  a dense image (see \cite[Proposition~0.4.6, p.~13]{arh1992} and \cite[Corollary~0.4.8, p.~13]{arh1992}). These results can be used in combination with Proposition~\ref{prop_caliber_dense_subset}(3) to establish the following fact:

\begin{proposition}\label{prop_Cp_condensation} If $X$ is condensed onto a Tychonoff space $Y$, then $\C(C_p(Y))\subseteq\C(C_p(X))$.

\end{proposition}

Later we will see that the inclusion $\C(C_p(Y))\subseteq\C(C_p(X))$ in Proposition~\ref{prop_Cp_condensation} can fail to be an equality after we have proven Proposition~\ref{prop_calibers_Cp_[0,kappa]}.

Now let $Y$ be a subspace of $X$. In \cite[Proposition~0.4.1, p.~11]{arh1992} it is shown that the restriction map from $C_p(X)$ onto $C_p(Y)$ is a continuous function with dense image. Since calibers are preserved by continuous maps, Proposition~\ref{prop_caliber_dense_subset} implies:

\begin{proposition}\label{prop_Cp_subspace} If $Y$ is a subspace of $X$, then $\C(C_p(X))\subseteq \C(C_p(Y))$.
\end{proposition}

A question that arises from Proposition~\ref{prop_Cp_subspace} is whether $\C(C_p(X))= \C(C_p(Y))$ for a subspace $Y$ with a specific property. Later, with the help of Proposition~\ref{prop_calibers_Cp_[0,kappa]}, we will see that this equality may fail even if $Y$ is $C$-embedded in $X$.

We recursively define the following spaces for each $n\in\mathbb{N}$: $$C_{p,1}(X):=C_p(X) \quad \text{and} \quad C_{p,n+1}(X) := C_p\left(C_{p,n}(X)\right).$$ It is shown in \cite[Proposition~0.5.5, p.~17]{arh1992} that $X$ is homeomorphic to a subspace of $C_{p,2}(X)$. On the other hand, it is known that $\C\left(X\right)\cap \UR =\C\left(C_{p,2}(X)\right)\cap \UR$ (see \cite[T.291, p.~319]{tkachuk2014}). The combination of these facts with Proposition~\ref{prop_Cp_subspace} implies:

\begin{corollary}\label{cor_Cp_calibers_iterated} For every $n\in\mathbb{N}$, $$\C\left(C_{p,n}(X)\right)\cap \UR \subseteq \C\left(C_{p,n+2}(X)\right)\subseteq \C\left(C_{p,n}(X)\right).$$

\end{corollary}

We will use the symbol $iw(X)$ to talk about the $i$-\textit{weight} of $X$ (see \cite{arh1992} and \cite{arhtka2008}). In \cite[Theorem~\textsc{I}.1.4, p.~26]{arh1992}, \cite[Theorem~\textsc{I}.1.5, p.~26]{arh1992} and \cite[Lemma~ 5.2.10, p.~299]{arhtka2008} it is proven that \begin{align}\label{fun_car_Cp(X)} d(X)=iw\left(C_p(X)\right) \quad \text{and} \quad d\left(C_p(X)\right)=iw(X)\leq nw(X).
\end{align}

With this in mind we can give some more results. For example, Proposition~\ref{prop_caliber_density} and~(\ref{fun_car_Cp(X)}) imply:

\begin{proposition}\label{prop_Cp_iw_Tychonoff} The relation $\{\kappa\in\CN : \cf(\kappa)>iw(X)\}\subseteq\C(C_p(X))$ holds.

\end{proposition}

An inductive argument combined with (\ref{fun_car_Cp(X)}) implies that, for every $n\in \mathbb{N}$, \begin{align}\label{fun_car_Cp(X)_iterado} d(X)=iw\left(C_{p,2n-1}(X)\right) \quad \text{and} \quad iw(X)=d\left(C_{p,2n-1}(X)\right).
\end{align} 	The next corollary follows from these equalities and Proposition~\ref{prop_Cp_iw_Tychonoff}.

\begin{corollary}\label{cor_prop_Cp_iw_Tychonoff} Let $n\in\mathbb{N}$. 
\begin{enumerate}
\item If $X$ has countable $i$-weight, then $C_{p,2n-1}(X)$ is strongly \v{S}anin.
\item When $X$ is separable, $C_{p,2n}(X)$ is strongly \v{S}anin.
\item Consequently, if $X$ is a separable space with countable $i$-weight, then $C_{p,n}(X)$ is strongly \v{S}anin.
\end{enumerate}

\end{corollary}

In particular,

\begin{corollary}\label{cor_Cp_compact_metrizable} If $X$ is metrizable and separable, then $C_{p,n}(X)$ is strongly \v{S}anin for every $n\in \mathbb{N}$.

\end{corollary}

In the realm of paracompactness (a natural weakening of metrizability), it is also possible to obtain some results. In Section~\ref{secc_Lindelof} we will calculate the calibers of $C_p\left(L(\lambda, \kappa)\right)$, where $L(\lambda,\kappa)$ is the $\lambda$-Lindelöf extension of a discrete space of cardinality $\kappa$. It turns out that these spaces are hereditarily paracompact and $\kappa^{+}\not \in \C(C_p(L(\lambda, \kappa)))$ (see Theorem~\ref{thm_Cp_calibers_lambda_Lindelof}). Thus, $C_p(L(\lambda, \kappa))$ is not a \v{S}anin space.

These spaces suggest a slightly more general construction. If $X$ is a topological space and $M$ is a subset of $X$, the {\it Bing-Hanner modification} is the space $X_M$ equipped with the topology generated by the base $\{U\cup A : U\in \tau_X \ \wedge \ A \subseteq M\}$. It is shown in \cite[Example~5.1.22, p.~306]{engelking1989} that if $X$ is a hereditarily paracompact Hausdorff space, then $X_M$ has the same characteristics. With the previous result it is possible to generalize in a certain way what happened with $C_p(L(\lambda, \kappa))$.

A {\it $\kappa$-sequence} in $X$ is a function $s: \kappa \to X$ usually denoted by $(x_\alpha)_{\alpha<\kappa}$ where $x_\alpha = s(\alpha)$ for every $\alpha<\kappa$. We say that a $\kappa$-sequence $(x_\alpha)_{\alpha<\kappa}$ {\it converges} to $x\in X$ if for every $U\in \tau_X(x)$ there is $\beta<\kappa$ with $\{x_\alpha : \alpha\geq \beta\} \subseteq U$.

\begin{proposition} Let $X$ be a hereditarily paracompact space, $x\in X$ and $\{x_\alpha : \alpha<\kappa\} \subseteq X\setminus \{x\}$ a $\kappa$-sequence convergent to $x$. If $M := \{x_\alpha : \alpha<\kappa\}$, then $X_M$ is hereditarily paracompact and $C_p(X_M)$ does not have caliber $\kappa$.

\end{proposition}

\begin{proof} Consider the family $\{V(\chi_\alpha, \{x_\alpha,x\},1/4) : \alpha<\kappa\}$, where $\chi_\alpha : X_M \to \mathbb{ R}$ is the function given by $$\chi_\alpha(y)= \begin{cases} 1, & \text{if} \ y = x_\alpha, \\
0, & \text{if} \ y\neq x_\alpha.
    \end{cases}$$ If we assume that $\kappa$ is a caliber for $C_p(X_M)$, then there are $J \in [\kappa]^{\kappa}$ and $f\in \bigcap\{V(\chi_\alpha, \{x_\alpha,x\},1/4) : \alpha\in J\}$. By convergence of the net $(x_\alpha)$, there exists $\beta<\kappa$ with $\{f(x_\alpha) : \alpha\geq \beta\} \subseteq (f(x)-1/ 2,f(x)+1/2)$. Hence, if $\alpha \in [\beta,\kappa) \cap J$, then $f\in V(\chi_\alpha, \{x_\alpha,x\},1/4)$ guarantees that $f(x) < 1/4$ and $3/4 < f(x_\alpha)$, which implies $|f(x)-f(x_\alpha)| \geq 1/2$, a contradiction to $f(x_\alpha) \in (f(x)-1/2,f(x)+1/2)$.
\end{proof}

On the other hand, by imposing some additional properties on a paracompact space $X$ it is possible to obtain some calibers for $C_p(X)$. For a point $x\in X$ we define the \textit{local $i$-weight} of $x$ in $X$ as the cardinal number $$liw(x,X)=\min\left\{iw(U) : U\in\tau_X(x)\right\}.$$ Hence, the \textit{local $i$-weight} of $X$ is $liw(X)=\sup\{liw(x,X) : x \in X\} $ (see \cite{bailey2005}).

In \cite[Theorem~4.23, p.~65]{bailey2005} it is shown that if $X$ is a paracompact Tychonoff space, then \begin{align}\label{iw_paracompacto} iw(X)=\max\left \{ \log(e(X)), liw(X)\right\}.
\end{align}

Additionally, it is known that if $X$ is a paracompact Hausdorff space, then $e(X)\leq c(X)$ (see \cite{hodel1984} and \cite[Lemma~4, p.~339 ]{sakai1983}). Finally, we say that $X$ is a \textit{locally ccc} space if for every $x\in X$ there exists $U\in \tau_X(x)$ with $c(U)=\omega$.

With these preliminaries we can obtain one more result:

\begin{proposition}\label{prop_Cp_paracompact} If $X$ is paracompact with countable local $i$-weight, then $$\left\{\kappa\in\CN : \cf(\kappa )>\log\left(c(X)\right)\right\}\subseteq \C(C_p(X)).$$
In particular, $C_p(X)$ is strongly \v{S}anin whenever $c(X)\leq \mathfrak{c}$.

\end{proposition}

\begin{proof} It follows from (\ref{iw_paracompacto}) that $iw(X)=\max\{\log(e(X)), liw(X)\} =\log(e(X))\leq \log(c(X))$. For this reason, Proposition~\ref{prop_Cp_iw_Tychonoff} yields $$\left\{\kappa\in\CN : \cf(\kappa)>\log\left(c(X)\right)\right\}\subseteq \left\{\kappa\in\CN : \cf(\kappa)>iw(X)\right\}\subseteq  \C(C_p(X)).$$ 
\end{proof}

Finally, it is sufficient that {\lq\lq}$c(X)\leq \mathfrak{c}${\rq\rq} in Proposition~\ref{prop_Cp_paracompact} to obtain that $C_p(X)$ is strongly \v{S}anin, but it is not necessary. For example, the Sorgenfrey line $\mathbb{S}$ is a submetrizable hereditarily paracompact space with countable $i$-weight; in particular, $C_p(\mathbb{S})$ is strongly \v{S}anin by Corollary~\ref{cor_prop_Cp_iw_Tychonoff}. Thus, if $\lambda$ is an infinite cardinal and $X$ is the product $\mathbb{S}\times D(\lambda)$, then $X$ is submetrizable, hereditarily paracompact (see \cite[Theorem~ 5.1.30, p.~308]{engelking1989}), has countable local $i$-weight, $c(X)=\lambda$ and $C_p(X)$ is strongly \v{S}anin due to Corollary~\ref{cor_Cp_sums}(3) that we will prove later on.

%%%%%%%%%%%%%%%%%%%%%%%%%%%%%%%%%%%%%%%%%%%%%%%%%%%%%%%%%%%%%%
\section{Calibers for $C_p(X)$ when $X$ is submetrizable}\label{secc_submetrizable}
%%%%%%%%%%%%%%%%%%%%%%%%%%%%%%%%%%%%%%%%%%%%%%%%%%%%%%%%%%%%%%

Tkachuk showed in \cite[Corollary~5, p.~270]{tkachuk1986} that if $X$ is a submetrizable space, then $C_p(X)$ is a \v{S}anin space. Our next goal is to improve Tkachuck's result by proving that, in fact, $C_p(X)$ is a strongly \v{S}anin space if $X$ is submetrizable.

If $\{X_\xi : \xi<\lambda\}$ is a family of topological spaces and $\{x_\alpha : \alpha<\kappa\}$ is a subset of $\prod_{\xi<\lambda} X_\xi$, the symbol $x_\alpha(\xi)$ represents the $\xi$-th coordinate of $x_\alpha$. For each $2\leq n < \omega$ let $$T_{X}^{n} := \left\{x\in X^n : \forall i<j<n(x(i) \neq x(j))\right\}.$$ On the other hand, $A\subseteq X$ is a {\it zero set} if there exists a continuous function $f: X \to [0,1]$ such that $A = f^{-1}\{0\}$. The family $Z(X)$ is the collection of zero subsets of $X$.

\begin{definition}\label{def_P(kappa)} Let $\kappa$ be a cardinal number and $2\leq n < \omega$. We say that {\it $X$ satisfies $P(\kappa,n)$} if for any $\{x_\alpha : \alpha<\kappa\} \subseteq T_{X}^{n}$ there exist $J \in [\kappa]^{\kappa}$ and a pairwise disjoint collection $\{C_i : i<n\} \subseteq Z(X)$ with $\{x_\alpha(i) : \alpha\in J\} \subseteq C_i$ for every $i<n$.

\end{definition}

With these preliminaries, Marsh proved in \cite[Theorem~60, p.~42]{marsh2004} the following result:

\begin{theorem}\label{thm_Marsh} If $\kappa$ is a cardinal number with $\cf(\kappa)>\omega$, then $C_p(X)$ has caliber $\kappa$ if and only if $X$ satisfies $P(\kappa,n)$ for every $2\leq n < \omega$.

\end{theorem}

The first we will show is that $P(\kappa,2)$ is sufficient to guarantee that $C_p(X)$ has caliber $\kappa$. Recall that $A\subseteq X$ and $B\subseteq X$ are {\it completely separated in $X$} if there exists a continuous function $f : X \to [0,1]$ such that $f[A]\subseteq \{0\}$ and $f[B] \subseteq \{1\}$.

\begin{proposition}\label{prop_Cp(X)_P(kappa)} If $\kappa$ is a cardinal number with $\cf(\kappa)>\omega$, then the following statements are equivalent.

\begin{enumerate}
\item $C_p(X)$ has caliber $\kappa$.
\item $X$ satisfies $P(\kappa,n)$ for each $2\leq n < \omega$.
\item $X$ satisfies $P(\kappa,2)$.
\item For every $\{(x_\alpha,y_\alpha) : \alpha<\kappa\} \subseteq T_X^2$ there is $J\in [\kappa]^{\kappa}$ such that $\{x_\alpha : \alpha\in J\}$ and $\{y_\alpha : \alpha\in J\}$ are completely separated in $X$.

\end{enumerate}

\end{proposition}

\begin{proof} The equivalence between (1) and (2) is Theorem~\ref{thm_Marsh}; the equivalence between (3) and (4) is Theorem (f), p. 37 in \cite{porwoo1988}, while (2) obviously implies (3). The proof that (3) implies (2) will be demostrated by induction on $n$. Suppose that for some $2\leq n < \omega$ the space $X$ satisfies $P(\kappa,2)$ and $P(\kappa,n)$. To show that $X$ also satisfies $P(\kappa,n+1)$ note that if $\{x_\alpha : \alpha<\kappa\}$ is a subset of $T_{X}^{n+1}$, then $P(\kappa,n)$ yields $I \in [\kappa]^{\kappa}$ and a pairwise disjoint collection $\{B_i : i<n\} \subseteq Z(X)$ with $\{x_\alpha(i) : \alpha\in I\} \subseteq B_i$ for any $i<n$.

\medskip

\noindent {\bf Claim.} There exist $\{J_i : i<n\} \subseteq [\kappa]^{\kappa}$ and $\{D(i,j) : i<n \ \wedge \ j<2\} \subseteq Z(X)$ such that $J_0 \subseteq I$ and if $i<n$, then $J_{i+1} \subseteq J_i$ (whenever $i+1<n$), $\{x_\alpha(n) : \alpha \in J_i\} \subseteq D(i,0)$, $\{x_\alpha(i) : \alpha \in J_i\} \subseteq D(i,1)$ and $D(i,0) \cap D(i,1) = \emptyset$.

\medskip

The proof of the statement above will be by finite recursion. To construct $J_0$ and $\{D(0,0),D(0,1)\}$ observe that $\{(x_\alpha(n),x_\alpha(0)) : \alpha \in I\} \subseteq T_{X}^{2}$ and, therefore, $P(\kappa,2)$ generates $J_0 \in [I]^{\kappa}$ and $\{D(0,0),D(0,1)\} \subseteq Z(X)$ with $\{x_\alpha(n) : \alpha \in J_0\} \subseteq D(0,0)$, $\{x_\alpha(0) : \alpha \in J_0\} \subseteq D(0,1)$ and $D(0,0)\cap D(0,1)= \emptyset$. Now, suppose that for some $i<n$ with $i+1<n$ we have constructed the collections $\{J_p : p\leq i\}$ and $\{D(p,j) : p\leq i \ \wedge \ j<2\}$ with the desired properties. Since $\{(x_\alpha(n),x_\alpha(i+1)) : \alpha \in J_i\} \subseteq T_{X}^{2}$, $P(\kappa,2)$ produces $J_{i+1} \in [J_i]^{\kappa}$ and $\{D(i+1,0),D(i+1,1)\} \subseteq Z(X)$ with $\{x_\alpha(n) : \alpha \in J_{i+1}\} \subseteq D(i+1,0)$, $\{x_\alpha(i+1) : \alpha \in J_{i+1}\} \subseteq D(i+1,1)$ and $D(i+1,0)\cap D(i+1,1)= \emptyset$. This completes the recursion.

To finish the proof of the proposition, define $C_n := \bigcap_{i<n} D(i,0)$ and $C_{i} := B_i \cap D(i,1)$ for every $i<n$. Clearly, $\{C_i : i<n+1\}$ is a pairwise disjoint family consisting of elements of $Z(X)$ such that $\{x_\alpha(i) : \alpha\in J_{n-1}\} \subseteq C_i$ for any $i<n+1$. In sum, $X$ satisfies $P(\kappa,n+1)$.
\end{proof}

Regarding the preservation of $P(\kappa,2)$ in topological products we have the following result:

\begin{proposition}\label{prop_P(kappa,2)_prod} If $\kappa$ and $\lambda$ are cardinal numbers with $\cf(\kappa)>\lambda\geq \omega$ and $\{X_\xi : \xi<\lambda\}$ is a family of spaces such that $C_p(X_\xi)$ has caliber $\kappa$ for every $\xi<\lambda$, then $C_p(\prod_{\xi<\lambda} X_\xi)$ has caliber $\kappa$.

\end{proposition}

\begin{proof} Set $X:= \prod_{\xi<\lambda} X_\xi$ and fix $\{(x_\alpha,y_\alpha) : \alpha<\kappa\} \subseteq T_{X}^{2}$. If $f : \kappa \to \lambda$ is the function determined by $f(\alpha) := \min\{\xi<\lambda : x_\alpha(\xi) \neq y_\alpha(\xi) \}$, then there exist $I\in [\kappa]^{\kappa}$ and $\xi<\lambda$ such that $f(\alpha) = \xi$ for any $\alpha \in I $. This implies that $\{(x_\alpha(\xi),y_\alpha(\xi)) : \alpha \in I\} \subseteq T_{X_\xi}^2$ and, therefore, there is $J\in [I]^{\kappa}$ such that $\{x_\alpha(\xi) : \alpha\in J\}$ and $\{y_\alpha(\xi) : \alpha\in J\}$ are completely separated in $X_\xi$. Hence, since the $\xi$-th projection is a continuous function, it turns out that $\{x_\alpha : \alpha\in J\}$ and $\{y_\alpha : \alpha\in J\}$ are completely separated in $X$.
\end{proof}

\begin{corollary}\label{cor_P(kappa,2)_prod} If $\{X_n : n<\omega\}$ is a family of spaces such that $C_p(X_n)$ is strongly \v{S}anin for every $n<\omega$, then $C_p(\prod_{n<\omega} X_n)$ is strongly \v{S}anin.

\end{corollary}

Now, let $\kappa$ be an infinite cardinal number and $J_\alpha := (0,1]\times \{\alpha\}$ for each $\alpha<\kappa$. Furthermore, let $J(\kappa ) := \{\langle 0,0\rangle\} \cup \bigcup_{\alpha<\kappa} J_\alpha$ and $d : J(\kappa) \to \mathbb{R}$ be the function determined by $$d\left( \langle x,\alpha \rangle,\langle y,\beta \rangle \right) := \begin{cases} \abs{x-y}, & \text{si} \ \alpha = \beta, \\ x+y, & \text{si} \ \alpha\neq \beta. \end{cases}$$ The space $J(\kappa)$ equipped with the metric $d$ is known as the {\it metrizable hedgehog of $\kappa$ spines}.

Finally, recall that if $A$ and $B$ are non-empty subsets of a metric space $(X,d)$, the {\it distance} between $A$ and $B$ is $$d(A,B) := \inf\{d(x,y) : x\in A \ \wedge \ y\in B\}.$$ Also, it is convenient to establish $d(A,\emptyset) := 1$ in order to simplify our next proof. 

In the proof of Theorem~\ref{thm_erizo_fuertemente_sanin} we use the following fact which will be proven in Corollary~\ref{cor_Cp_sums}: if $\lambda$ is a cardinal number and $X$ is a topological space, then $\C(C_p(X)) = \C(C_p(X\times D(\lambda)))$.

\begin{theorem}\label{thm_erizo_fuertemente_sanin} $C_p(J(\kappa))$ is strongly \v{S}anin for every infinite cardinal $\kappa$.

\end{theorem}

\begin{proof} Let $\lambda$ be a cardinal with $\cf(\lambda)>\omega$ and $\{V(f_\beta,A_\beta,\varepsilon) : \beta<\lambda\}$ be an elementary family in $C_p(J(\kappa))$ with $\langle 0,0 \rangle\in \bigcap_{\beta<\lambda} A_\beta$. Since $\mathbb{R}$ has caliber $\lambda$ and $\{B_\varepsilon(f_\beta(\langle 0,0 \rangle)) : \beta<\lambda\}\subseteq \tau_\mathbb{ R}^+$, there exist $I\in [\lambda]^{\lambda}$ and $r\in \bigcap\{B_\varepsilon(f_\beta(\langle 0,0 \rangle)) : \beta \in I\}$. Then, since each $f_\beta$ is continuous and $f_\beta(\langle 0,0 \rangle) \in B_{\varepsilon}(r)$, there exists $0<\delta_\beta<1$ with $f_\beta[\{\langle 0,0 \rangle\}\cup ((0,\delta_\beta]\times \kappa)] \subseteq B_{\varepsilon}(r)$. Thus, as $\mathbb{R}$ has caliber $\lambda$ and $\{(0,\delta_\beta) : \beta \in I\} \subseteq \tau_\mathbb{R}^+$, there are $J\in [I]^{\lambda}$ and $\delta \in \bigcap\{(0,\delta_\beta) : \beta \in J\}$.

Let $p<\omega$ satisfy $0<\delta-2^{-p}<\delta+2^{-p}<1$, and set $B := \{\delta\} \times \kappa$ and $B_\beta := A_\beta \setminus (\{\langle 0,0 \rangle\}\cup ((0,\delta]\times \kappa))$ for each $\beta \in J$. Since $B$ is closed, $B_\beta$ is finite and $B\cap B_\beta = \emptyset$, the function $h : J \to \omega$ given by $h(\beta) := \min \{p\leq n<\omega : d(B,B_\beta) > 2^{-n}\}$ is well defined. Moreover, there exists $p\leq m<\omega$ such that if $K := h^{-1}\{m\}$, then $K\in [J]^{\lambda}$ and $d(B,B_\beta) > 2^{-m}$ for any $\beta \in K$.

Let $\eta := \delta+2^{-m}$ and $X := [\eta,1]\times \kappa$. Since the identity function $X \to [\eta,1] \times D(\kappa)$ is a homeomorphism, it follows that $\C(C_p(X)) = \C(C_p([\eta,1] \times D(\kappa))) = \C(C_p([\eta,1])) = \UC$ (see Proposition~\ref{prop_Cp_iw_Tychonoff} and Corollary~\ref{cor_Cp_sums}). For this reason, since $\{V(\restr{f_\beta}{X}, B_\beta,\varepsilon) : \beta \in K\}$ is a family of open sets in $C_p(X)$, there exist $L\in [K]^{\lambda}$ and a continuous function $g: X\to \mathbb{R}$ with $g\in \bigcap\{V(\restr{f_\beta}{X}, \beta,\varepsilon) : \beta \in L\}$.

Now, for each $\alpha<\kappa$ let $g_\alpha : [\delta,\eta]\times\{\alpha\} \to \mathbb{R}$ be the function given by $$g_\alpha\left(\langle x,\alpha\rangle\right) := \left(\frac{g\left(\langle \eta,\alpha\rangle\right)-r}{\eta-\delta}\right) (x-\delta)+r.$$ Furthermore, observe that the function $f : J(\kappa) \to \mathbb{R}$ determined by $$f\left(\langle x,\alpha\rangle\right) :=  \begin{cases} r, & \ \text{si} \ x\in [0,\delta], \\
g_\alpha\left(\langle x,\alpha\rangle\right), & \ \text{si} \ x\in [\delta,\eta], \\
g\left(\langle x,\alpha\rangle\right), & \ \text{si} \ x\in [\eta,1],
\end{cases}$$ is continuous.

To finish our proof, notice that if $\beta \in L$, $\langle x,\alpha \rangle \in A_\beta$ and $z:= \langle x,\alpha \rangle$ two cases arise: if $0\leq x \leq \delta$, then $f_\beta(z) \in B_\varepsilon(r)$; that is, $|f(z) - f_\beta(z)| = |r-f_\beta(z)|<\varepsilon$. If $x>\delta$, then $z\in B_\beta$, which implies that $$x-\delta = |x-\delta| = d\left(\langle \delta,\alpha\rangle,\langle x,\alpha\rangle\right) = d\left(\langle \delta,\alpha\rangle,z\right) \geq d(B,B_\beta) > 2^{- m};$$ i.e., $x>\delta+2^{-m} = \eta$. Therefore, $|f(z) - f_\beta(z)| = |g(z)-f_\beta(z)|<\varepsilon$. In sum, $f\in \bigcap\{V(f_\beta,A_\beta,\varepsilon) : \beta\in L\}$.
\end{proof}

A consequence of Corollary~\ref{cor_P(kappa,2)_prod} and Theorem~\ref{thm_erizo_fuertemente_sanin} is:

\begin{theorem}\label{thm_potencia_erizo_fuertemente_Sanin} $C_p(J(\kappa)^\omega)$ is strongly \v{S}anin for every infinite cardinal $\kappa$. 

\end{theorem}

Hence,

\begin{theorem}\label{thm_Cp(X)_metrizable_f_sanin} If $X$ is a metrizable space, then $C_p(X)$ is strongly \v{S}anin.

\end{theorem}

\begin{proof} By Kowalsky's Theorem (see \cite[Theorem~4.4.9, p.~282]{engelking1989}), there exists an infinite cardinal $\kappa$ such that $X$ is a subspace of the countable power $J(\kappa)^\omega$. Then, since $C_p(J(\kappa)^\omega)$ is strongly \v{S}anin by Theorem~\ref{thm_potencia_erizo_fuertemente_Sanin}, Proposition~\ref{prop_Cp_subspace} guarantees that $C_p(X)$ is strongly \v{S}anin.
\end{proof}

Finally, Proposition~\ref{prop_Cp_condensation} and Theorem~\ref{thm_Cp(X)_metrizable_f_sanin} produce the following result:

\begin{theorem}\label{thm_Cp(X)_submetrizable_f_sanin} If $X$ is a submetrizable space, then $C_p(X)$ is strongly \v{S}anin.

\end{theorem}

%%%%%%%%%%%%%%%%%%%%%%%%%%%%%%%%%%%%%%%%%%%%%%%%%%%%%%%%%%%%%%
\section{Calibers for $C_p(X)$ when $X$ is an interval of ordinals}\label{secc_ordinales}
%%%%%%%%%%%%%%%%%%%%%%%%%%%%%%%%%%%%%%%%%%%%%%%%%%%%%%%%%%%%%%

If $\kappa$ is a cardinal number with $\cf(\kappa)>\omega$, the functions in $C\left([0,\kappa]\right)$ and $C\left( [0,\kappa)\right)$ are constant in an interval of the form $[\alpha,\kappa)$ for some $\alpha<\kappa$ (see \cite[Theorem~(q)(5), p.~130]{porwoo1988}). This fact allows us to establish the following result:

\begin{lemma}\label{lema_Cp_[0,kappa]_no_caliber_kappa} If $\kappa$ is an infinite cardinal, then $\kappa$ is not a caliber for $C_p\left([0,\kappa]\right)$ nor for $C_p\left([0,\kappa)\right )$.

\end{lemma}

\begin{proof} Since both arguments are similar, we focus on $C_p\left([0,\kappa]\right)$. If $\cf(\kappa)=\omega$, then Proposition~\ref{cor_cellularity_cofinality} confirms that $C_p\left([0,\kappa]\right)$ does not have caliber $\kappa$. When $\cf(\kappa)>\omega$, for each $\alpha\in S := \{\beta+1 : \beta<\kappa\}$ we define $\chi_\alpha:[0,\kappa ] \to\mathbb{R}$ by $$\chi_\alpha(\beta)= \begin{cases} 1, & \text{if} \ \beta=\alpha, \\
0, & \text{if} \ \beta\neq\alpha.
   \end{cases}$$ Consider the collection $\left\{V(\chi_\alpha, \{\alpha,\alpha+1\}, 1/2) : \alpha\in S\right\}$ and let $J\in [S]^{\kappa}$. If there exists $f\in \bigcap\{V(\chi_\alpha, \{\alpha,\alpha+1\}, 1/2) : \alpha\in J\}$, then there are $\beta <\kappa$ and $r\in\mathbb{R}$ with $f\left[[\beta,\kappa]\right] \subseteq \{r\}$. Finally, if $\alpha\in [\beta,\kappa]\cap J$, then $f\in V(\chi_\alpha, \{\alpha,\alpha+1\}, 1/2 )$ implies that $1/2<r$ and $r<1/2$ simultaneously, an absurdity.
\end{proof}

With the previous lemma we can prove the following proposition:

\begin{proposition}\label{prop_calibers_Cp_[0,kappa]} If $\kappa$ is an infinite cardinal, then $$\C\left(C_p\left([0,\kappa]\right)\right)=\{\lambda\in\CN : \cf(\lambda)>\kappa\}=\C\left(C_p\left([0,\kappa)\right)\right).$$
\end{proposition}

\begin{proof} Let $X:= C_p\left([0,\kappa]\right)$ and $Y:= C_p\left([0,\kappa)\right)$. First, $nw\left([0,\kappa]\right)=\kappa=nw\left([0,\kappa)\right)$ and $nw\geq iw$ imply that $\{\lambda \in\CN : \cf(\lambda)>\kappa\}\subseteq \C(X)\cap \C(Y)$ (see Proposition~\ref{prop_Cp_iw_Tychonoff}). On the other hand, if $\lambda\in \CN$ satisfies $\cf(\lambda)\leq\kappa$, then since $C_p\left([0,\cf(\lambda) )\right)$ does not have $\cf(\lambda)$ as a caliber (see Lemma~\ref{lema_Cp_[0,kappa]_no_caliber_kappa}), $X$ and $Y$ also do not have caliber $\cf(\lambda)$ (see Proposition~\ref{prop_Cp_subspace}). Hence, Proposition~\ref{prop_caliber_cofinality} ensures that $\lambda \not \in \C(X)\cup \C(Y)$.
\end{proof}

In particular, if $\kappa$ and $\lambda$ are cardinals with $\omega \leq \lambda<\kappa$, then $[0,\lambda]$ is $C$-embedded in $[0,\kappa]$, $\lambda^{+}$ is a caliber for $C_p\left([0,\lambda]\right)$ and is not a caliber for $C_p\left([0,\kappa]\right)$. This shows that the other inclusion in Proposition~\ref{prop_Cp_subspace} is not true in general, even if the subspace is $C$-embedded.

On the other hand, since the discrete space $D(\kappa)$ is condensed onto $[0,\kappa]$ and $\kappa \in \C\left(C_p\left(D(\kappa)\right)\right) \setminus \C\left(C_p\left([0,\kappa]\right)\right)$
by Propositions~\ref{prop_Cp_calibers_discrete} and~\ref{prop_calibers_Cp_[0,kappa]}, it follows that the other inclusion in Proposition~\ref{prop_Cp_condensation} is also not true in general.

We can further strengthen Proposition~\ref{prop_calibers_Cp_[0,kappa]} for any initial segment of ordinals.

\begin{proposition}\label{prop_calibers_Cp_[0,alpha]} If $\alpha$ is an infinite ordinal, then $$\C\left(C_p\left([0,\alpha]\right)\right) = \left\{\lambda\in \CN : \cf(\lambda)>\abs{\alpha}\right\} = \C\left(C_p\left([0,\alpha)\right)\right).$$

\end{proposition}

\begin{proof} When $\abs{\alpha} = \alpha$ the result is Proposition~\ref{prop_calibers_Cp_[0,kappa]}. If $\abs{\alpha} < \alpha$ we observe that, if $X := C_p\left(\left[0,\alpha\right]\right)$ and $Y := C_p\left(\left[ 0,\alpha\right)\right)$, then the relations $$iw\left(\left[0,\alpha\right]\right) \leq nw\left(\left[0,\alpha\right] \right) \leq |\alpha| \quad \text{and} \quad iw\left(\left[0,\alpha\right)\right) \leq nw\left(\left[0,\alpha\right)\right) \leq |\alpha |$$ and Proposition~\ref{prop_Cp_iw_Tychonoff} imply that $\left\{\lambda\in \CN : \cf(\lambda)>\abs{\alpha}\right\}\subseteq \C(X) \cap \C(Y)$. Finally, since $\left[0,\abs{\alpha}\right]$ is a subspace of $[0,\alpha)$ and $[0,\alpha]$, Propositions~\ref{prop_Cp_subspace} and~\ref{prop_calibers_Cp_[0,kappa]} confirm that $\C(X)\cup \C(Y) \subseteq \left\{\lambda\in \CN : \cf(\lambda )>\abs{\alpha}\right\}$.
\end{proof}

%%%%%%%%%%%%%%%%%%%%%%%%%%%%%%%%%%%%%%%%%%%%%%%%%%%%%%%%%%%%%%
\section{Calibers for $C_p(X)$ when $X$ is the one-point $\lambda$-Lindelöf extension}\label{secc_Lindelof}
%%%%%%%%%%%%%%%%%%%%%%%%%%%%%%%%%%%%%%%%%%%%%%%%%%%%%%%%%%%%%%

Let us now concentrate on the one-point $\lambda$-Lindelöf extension of an infinite discrete space. Recall that for an infinite cardinal $\kappa$ the discrete space $D(\kappa)$ is the set $\kappa$ equipped with the discrete topology. Now, if $\lambda$ is a cardinal number with $\omega \leq \lambda \leq \kappa$, the {\it one-point $\lambda$-Lindelöf extension} of the space $D(\kappa)$ is the set $L (\lambda,\kappa) := \kappa\cup \{\kappa\}$ endowed with the topology $$\tau_{L(\lambda,\kappa)} = \tau_{D(\kappa)} \cup \left\{U\subseteq L(\lambda,\kappa) : \kappa \in U \ \wedge \ \abs{L(\lambda,\kappa)\setminus U}<\lambda\right\}.$$ Observe that $L(\omega,\kappa)$ is the one-point compactification of $D(\kappa)$ and $L(\omega_1,\kappa)$ is the one-point Lindelöf extension of $D(\kappa)$.

\begin{lemma}\label{lemma_Cp_L(lambda,kappa)_no_caliber_mu} If $\kappa$, $\lambda$ and $\mu$ are cardinals with $\omega \leq \lambda\leq \mu\leq \kappa$, then $C_p\left(L(\lambda,\kappa)\right )$ does not have caliber $\mu$.

\end{lemma}

\begin{proof} For each $\alpha<\mu$ denote by $\chi_\alpha: L(\lambda,\kappa)\to\mathbb{R}$ the characteristic function of $\{\alpha\}$. Furthermore, for all $i<2$ define $F_i := [i-1/4, i+1/4]$. Consider $\{V(\chi_\alpha, \{\alpha,\alpha+1\}, 1/4) : \alpha<\mu\}$, take $J\in [\mu]^{\mu}$ and suppose, in search of a contradiction, that there exists $f\in \bigcap\left\{V(\chi_\alpha, \{\alpha,\alpha+1\}, 1/4) : \alpha \in J\right\}$. Now, since $F_0$ and $F_1$ are disjoint subsets of $\mathbb{R}$, there exists $i<2$ such that $f(\kappa)\not\in F_i$; consequently, there exists $\varepsilon>0$ with $B_\varepsilon(f(\kappa)) \cap F_i=\emptyset$. Let $U\in \tau_{L(\lambda,\kappa)}(\kappa)$ be such that $f[U] \subseteq B_\varepsilon(f(\kappa))$.

Observe that if $\abs{J\cap U} < \mu$, then $\abs{J\setminus U} = \mu$, which implies that $\mu = \abs{J\setminus U} \leq \abs{L(\lambda,\kappa) \setminus U}<\lambda\leq \mu$; a contradiction. Thus, $\abs{J\cap U} = \mu$ and therefore, the set $\{\alpha+1 : \alpha \in J\cap U\}$ has cardinality $\mu$. Let $ \alpha \in J\cap U$ be such that $\alpha+1\in U$. Hence, since $f\in V(\chi_\alpha, \{\alpha,\alpha+1\}, 1/4)$, it turns out that $f(\alpha) \in B_\varepsilon(f( \kappa)) \cap F_1$ and $f(\alpha+1) \in B_\varepsilon(f(\kappa)) \cap F_0$; a contradiction to $B_\varepsilon(f(\kappa)) \cap F_i=\emptyset$.
\end{proof}

A combination of Proposition~\ref{prop_caliber_cofinality} and Lemma~\ref{lemma_Cp_L(lambda,kappa)_no_caliber_mu} produces the following corollary.

\begin{corollary}\label{cor_Cp_L(lambda,kappa)_no_caliber_mu} If $\kappa$, $\lambda$ and $\mu$ are cardinals with $\omega\leq \lambda\leq \cf(\mu)\leq \kappa$, then $C_p\left(L(\lambda,\kappa)\right)$ does not have caliber $\mu$.

\end{corollary}

In our next result we will use the equivalence of Lemma~\ref{lemma_elementary_family}.

\begin{lemma}\label{lemma_Cp_L(lambda,kappa)_caliber_mu_1} If $\kappa$, $\lambda$ and $\mu$ are cardinals with $\mu$ regular and $\omega<\mu<\lambda\leq \kappa$, then $C_p\left(L(\lambda,\kappa)\right)$ has caliber $\mu$.

\end{lemma}

\begin{proof} Let $\{V(f_\alpha,A_\alpha,\varepsilon) : \alpha<\mu\}$ be an elementary family. By Theorem~\ref{thm_regular_delta_system}, there exist $I\in [\mu]^{\mu}$ and a set $A$ with $A = A_\alpha \cap A_\beta$, whenever $\alpha ,\beta \in I$ are different. For the rest of the argument we will assume, without loss of generality, that $A\neq \emptyset$ and that $A_\alpha \setminus A \neq \emptyset$ for every $\alpha\in I$. If this is not the case, since $|\bigcup_{\alpha<\mu} A_\alpha| \leq \mu<\kappa$, there is a set $\{\xi_\alpha : \alpha<\mu + 1\} \subseteq \kappa \setminus \bigcup_{\alpha<\mu} A_\alpha$ enumerated without repetitions. Now, we define $B_\alpha := A_\alpha \cup\{\xi_\alpha,\xi_\mu\}$, and if $I\in [\mu]^ {\mu}$ and $B$ satisfy $B = B_\alpha \cap B_\beta$, whenever $\alpha ,\beta \in I$ are different, then $\xi_\mu\in B$ and $ \xi_\alpha\in B_\alpha \setminus B$ for any $\alpha\in I$. Furthermore, for every $\alpha\in I$ the relation $f\in V(f_\alpha,B_\alpha,\varepsilon)$ implies that $f\in V(f_\alpha,A_\alpha,\varepsilon)$.

That said, let $m, n\in \omega$ be such that $A = \{\xi_i : i<m\}$ and $A_\alpha \setminus A = \{\xi(\alpha,i ) : i<n\}$ for each $\alpha \in I$.

\medskip

\noindent {\bf Case 1.} $\kappa \not\in \bigcup\{A_\alpha : \alpha\in I\}$.

\medskip

Since the following reasoning will be required several times, we present all the details of its proof.

\medskip

\noindent {\bf Claim.} There are collections $\{E_i : i<m\}$, $\{x_i : i<m\}$, $\{F_j : j<n\}$ and $\{y_j : j<n\}$ such that the following conditions hold for every $i<m$ and $j<n$:

\begin{enumerate}
\item $E_0 \in [I]^{\mu}$,
\item if $i+1<m$, then $E_{i+1} \in [E_i]^{\mu}$,
\item $x_i \in \bigcap \{B_\varepsilon(f_\alpha(\xi_i)) : \alpha \in E_i\}$,
\item $F_0 \in [E_{m-1}]^{\mu}$,
\item if $j+1<n$, then $F_{j+1} \in [F_j]^{\mu}$, and
\item $y_j \in \bigcap \{(B_\varepsilon(f_\alpha(\xi(\alpha,j))) : \alpha \in F_j\}$.
\end{enumerate}

\medskip

We obtain items (1), (2), (3) above by finite recursion on $m$. Indeed, since $\mathbb{R}$ has caliber $\mu$ by Proposition~\ref{prop_Cp_precalibers}(1), there exist $E_0 \in [I]^{\mu}$ and $x_0 \in \bigcap \{B_\varepsilon(f_\alpha(\xi_0)) : \alpha \in E_0\}$. Then, if for some $i$ with $i+1<m$ we have constructed families $\{E_j : j\leq i\}$ and $\{x_j : j\leq i\}$ with the desired properties, then the fact that $\mathbb{R}$ has caliber $\mu$ implies the existence of $E_{i+1} \in [E_i]^{\mu}$ and $x_{i+1} \in \bigcap \{B_\varepsilon(f_\alpha(\xi_{i+1})) : \alpha \in E_{i+1}\}$.

We arrive to the remaining three points by finite recursion on $n$. To begin with, the fact that $\mathbb{R}$ has caliber $\mu$ implies that there are $F_0 \in [E_{m-1}]^{\mu}$ and $y_0 \in \bigcap \{ B_\varepsilon(f_\alpha(\xi(\alpha,0))) : \alpha \in F_0\}$. What follows is to assume that for some $j$ with $j+1<n$ we have constructed families $\{F_k : k\leq j\}$ and $\{y_k : k\leq j\}$ with the required characteristics. Finally, since $\mathbb{R}$ has caliber $\mu$, it turns out that there exist $F_{j+1} \in [F_j]^{\mu}$ and $y_{j+1} \in \{ B_\varepsilon(f_\alpha(\xi(\alpha,j+1))) : \alpha \in F_{j+1}\}$.

Now, let $f: L(\lambda,\kappa) \to \mathbb{R}$ be the function given by $$f(\xi) := \begin{cases} 0, & \text{if} \ \xi \in L(\lambda,\kappa) \setminus \bigcup\{A_\alpha : \alpha\in F_{n-1}\}, \\
x_i, & \text{if} \ \exists i<m\left(\xi = \xi_i\right), \\
y_j, & \text{if} \ \exists\alpha\in F_{n-1} \exists j<n\left(\xi = \xi(\alpha,j)\right).
   \end{cases}$$ Observe that the relations $\abs{\bigcup\{A_\alpha : \alpha\in F_{n-1}\}} \leq \mu < \lambda$ and $\kappa \not \in \bigcup\{A_\alpha : \alpha\in F_{n-1}\}$ imply that $f$ is a continuous function. Furthermore, $f$ is an element of $\bigcap\{V(f_\alpha,A_\alpha,\varepsilon) : \alpha \in F_{n-1}\}$.

\medskip

\noindent {\bf Case 2.} $\kappa \in \bigcup\{A_\alpha : \alpha\in I\}$.

\medskip

Let $I^{*} := \{\alpha\in I : \kappa \in A_\alpha\}$. When $\abs{I\setminus I^{*}} = \mu$, then $\kappa \not\in \bigcup\{A_\alpha : \alpha\in I\setminus I^{*}\}$. Therefore, a similar reasoning to that performed for the Claim of Case 1 produces $J\in [I\setminus I^{*}]^{\mu}$ and $f\in C_p(L(\lambda,\kappa ))$ with $f\in \bigcap\{V(f_\alpha,A_\alpha,\varepsilon) : \alpha \in J\}$. Now suppose that $\abs{I\setminus I^{*}} <  \mu$. Since $\kappa \in \bigcap\{A_\alpha : \alpha\in I^{*}\}$, we can assume without loss of generality that $\xi_0 = \kappa$. Then, an argument similar to that made in Case 1 generates collections $\{E_i : i<m\}$, $\{x_i : i<m\}$, $\{F_j : j<n\}$ and $\{y_j : j<n\}$ such that the following conditions hold for every $i<m$ and $j<n$:

\begin{enumerate}
\item $E_0 \in [I^{*}]^{\mu}$ and $F_0 \in [E_{m-1}]^{\mu}$,
\item if $i+1<m$ and $j+1<n$, then $E_{i+1} \in [E_i]^{\mu}$ and $F_{j+1} \in [F_j]^{\mu}$,
\item $x_i \in \bigcap \{B_\varepsilon(f_\alpha(\xi_i)) : \alpha \in E_i\}$,
\item $y_j \in \bigcap \{(B_\varepsilon(f_\alpha(\xi(\alpha,j))) : \alpha \in F_j\}$.
\end{enumerate}

Finally, consider the function $f: L(\lambda,\kappa) \to \mathbb{R}$ defined by $$f(\xi) := \begin{cases} x_0, & \text{if} \ \xi \in L(\lambda,\kappa) \setminus \bigcup\{A_\alpha : \alpha\in F_{n-1}\}, \\
x_i, & \text{if} \ \exists i<m\left(\xi = \xi_i\right), \\
y_j, & \text{if} \ \exists\alpha\in F_{n-1} \exists j<n\left(\xi = \xi(\alpha,j)\right).
   \end{cases}$$ It remains to notice that $\abs{\bigcup\{A_\alpha\setminus \{\kappa\} : \alpha\in F_{n-1}\}} \leq \mu < \lambda$, $\kappa \not\in \bigcup\{A_\alpha\setminus\{\kappa\} : \alpha\in F_{n-1}\}$ and $f[L(\lambda,\kappa) \setminus \bigcup\{A_\alpha\setminus \{\kappa\} : \alpha\in F_{n-1}\}] \subseteq \{x_0\}$ imply that $f$ is a continuous function and $f\in \bigcap\{V(f_\alpha,A_\alpha,\varepsilon) : \alpha \in F_{n-1}\}$.
\end{proof}

\begin{corollary}\label{cor_Cp_L(lambda,kappa)_pi_peso} If $\kappa$, $\lambda$ and $\mu$ are cardinals with $\omega<\cf(\mu)<\lambda\leq \kappa<\mu$, then $C_p\left(L (\lambda,\kappa)\right)$ has caliber $\mu$.
\end{corollary}

\begin{proof} Since $\kappa$ is a dense, discrete and open subspace of $L(\lambda,\kappa)$, $\{\{\alpha\} : \alpha<\kappa\}$ is a $\pi$-base for $L(\lambda,\kappa)$ of size $\kappa$; consequently, $\pi w\left(L(\lambda,\kappa)\right) \leq \kappa<\mu$. Then, since $\cf(\mu)$ is a caliber for $C_p(L(\lambda,\kappa))$ by Lemma~\ref{lemma_Cp_L(lambda,kappa)_caliber_mu_1}, it follows from Proposition~\ref{prop_caliber_pi_weight} that $\mu$ is a caliber for $C_p(L(\lambda,\kappa))$.
\end{proof}

Our next result is the version for singular cardinals of Lemma~\ref{lemma_Cp_L(lambda,kappa)_caliber_mu_1}.

\begin{lemma}\label{lemma_Cp_L(lambda,kappa)_caliber_mu_2} If $\kappa$, $\lambda$ and $\mu$ are cardinals with $\omega<\cf(\mu)<\mu<\lambda\leq \kappa$, then $C_p\left(L (\lambda,\kappa)\right)$ has caliber $\mu$.
\end{lemma}

\begin{proof} Let $\{V(g_\gamma,B_\gamma,\varepsilon) : \gamma<\mu\}$ be an elementary family. Furthermore, let $\{\mu_\alpha : \alpha<\cf(\mu)\}$ be a collection of cardinals with the following properties:
\begin{enumerate}
\item $\mu_0 = 0$ and $\sup\{\mu_\alpha : \alpha<\cf(\mu)\}=\mu$, and
\item for each $0 < \alpha<\cf(\mu)$, $\mu_\alpha$ is a regular and uncountable cardinal such that $\max\left\{\cf(\mu), \sup\{\mu_\beta : \beta<\alpha\}\right\}<\mu_\alpha$.
\end{enumerate}

Now, let $I_\alpha := [\mu_\alpha,\mu_{\alpha+1})$ for every $\alpha<\cf(\mu)$. For each $\gamma \in \bigcup_{\alpha<\cf(\mu)} I_\alpha$ let $f(\alpha,\gamma) := g_\gamma$ and $A(\alpha,\gamma) := B_\gamma$, where $\alpha<\cf( \mu)$ is the only ordinal with $\gamma \in I_\alpha$. By Theorem~\ref{thm_singular_delta_system}, there exist $I\in \left[\cf(\mu)\right]^{\cf(\mu)}$, $\{J_\alpha : \alpha\in I\}$, $\{A_\alpha : \alpha\in I\}$ and $A$ that satisfy the following conditions for any different $\alpha,\beta\in I$:
\begin{enumerate}
\item $J_\alpha \in [I_\alpha]^{\mu_{\alpha+1}}$,
\item $A_\alpha \cap A_\beta = A$,
\item if $\gamma,\delta \in J_\alpha$ are different, then $A(\alpha, \gamma) \cap A(\alpha, \delta) = A_\alpha$,
\item when $\gamma\in J_\alpha$ and $\delta\in J_\beta$, then $A(\alpha, \gamma) \cap A(\beta, \delta) = A$, and
\item $|A_\alpha|=|A_\beta|$.
\end{enumerate}

For the rest of the proof we assume without losing generality that $A\neq \emptyset$, that $A_\alpha \setminus A \neq \emptyset$ for each $\alpha\in I$, and that $A(\alpha,\gamma)\setminus A_\alpha \neq \emptyset$ for every $\alpha \in I$ and $\gamma \in J_\alpha$. Let $m, n, p\in \omega$ be such that $A = \{\xi_i : i<m\}$, $A_\alpha \setminus A = \{\xi(\alpha,j) : j<n\}$ for each $\alpha \in I$, and $A(\alpha,\gamma)\setminus A_\alpha = \{\xi(\alpha,\gamma,k) : k<p\}$ for every $\alpha \in I$ and $\gamma \in J_\alpha$.

\medskip

\noindent {\bf Case 1.} $\kappa \not \in \bigcup\{A(\alpha, \gamma) : \alpha \in I \ \wedge \ \gamma \in J_\alpha\}$.

\medskip

Below we present three statements without proof since each is proven with a procedure similar to that set forth in the Claim of Lemma~\ref{lemma_Cp_L(lambda,kappa)_caliber_mu_1}.

\begin{enumerate}
\item For any $\alpha \in I$ there exist collections $\{E(\alpha,i) : i<m\}$ and $\{x(\alpha,i) : i<m\}$ such that the following conditions are satisfied for every $i<m$:

\begin{enumerate}
\item $E(\alpha,0) \in [J_\alpha]^{\mu_{\alpha+1}}$,
\item if $i+1<m$, then $E(\alpha,i+1) \in [E(\alpha,i)]^{\mu_{\alpha+1}}$, and
\item $x(\alpha,i) \in \bigcap \{B_{\varepsilon/2}(f(\alpha,\gamma)(\xi_i)) : \gamma \in E(\alpha,i)\}$.
\end{enumerate}

\item There are collections $\{K_i : i<m\}$ and $\{x_i : i<m\}$ such that the following conditions are satisfied for every $i<m$:

\begin{enumerate}
\item $K_0 \in [I]^{\cf(\mu)}$,
\item if $i+1<m$, then $K_{i+1} \in [K_i]^{\cf(\mu)}$, and
\item $x_i \in \bigcap \{B_{\varepsilon/2}(x(\alpha,i)) : \alpha \in K_i\}$.
\end{enumerate}

\item For any $\alpha \in K_{m-1}$ there exist collections $\{F(\alpha,j) : j<n\}$, $\{y(\alpha,j) : j<n \}$, $\{G(\alpha,k) : k<p\}$ and $\{z(\alpha,k) : k<p\}$ such that the following conditions are satisfied for every $i<m$, $j<n$ and $k<p$:

\begin{enumerate}
\item $F(\alpha,0) \in [E(\alpha,m-1)]^{\mu _{\alpha+1}}$ and $G(\alpha,0) \in [F(\alpha,n-1)]^{\mu_{\alpha+1}}$,
\item if $j+1<n$ and $k+1<p$, then $F(\alpha,j+1) \in [F(\alpha,j)]^{\mu_{\alpha+1 }}$ and $G(\alpha,k+1) \in [G(\alpha,k)]^{\mu_{\alpha+1}}$,
\item $y(\alpha,j) \in \bigcap \{B_\varepsilon(f(\alpha,\gamma)(\xi(\alpha,j))) : \gamma \in F(\alpha,j )\}$, and
\item $z(\alpha,k) \in \bigcap \{B_\varepsilon(f(\alpha,\gamma)(\xi(\alpha,\gamma,k))) : \gamma \in G(\alpha,k)\}$.

\end{enumerate}

\end{enumerate}

Now let $f: L(\lambda,\kappa) \to \mathbb{R}$ be the function given by $$f(\xi) := \begin{cases} 0, & \text{if} \ \xi \in L(\lambda,\kappa) \setminus \bigcup\{A(\alpha,\gamma) : \alpha \in K_{m-1} \ \wedge \ \gamma\in G(\alpha,p-1)\}, \\
x_i, & \text{if} \ \exists i<m\left(\xi = \xi_i\right), \\
y(\alpha,j), & \text{if} \ \exists\alpha\in K_{m-1} \exists j<n\left(\xi = \xi(\alpha,j)\right), \\
z(\alpha,k), & \text{if} \ \exists\alpha\in K_{m-1} \exists \gamma \in G(\alpha,p-1) \exists k<p\left(\xi = \xi(\alpha,\gamma,k)\right).
\end{cases}$$ Observe that the relations $\abs{\bigcup\{A(\alpha,\gamma) : \alpha \in K_{m-1} \ \wedge \ \gamma\in G(\alpha,p-1)\} } \leq \mu < \lambda$ and $\kappa \not\in \bigcup\{A(\alpha,\gamma) : \alpha \in K_{m-1} \ \wedge \ \gamma\in G( \alpha,p-1)\}$ imply that $f$ is a continuous function. Finally, notice that $G := \bigcup\{G(\alpha, p-1) : \alpha \in K_{m-1}\}$ satisfies $\abs{G} = \mu$ and, furthermore, $f$ is an element of $\bigcap\{V(g_\gamma,B_\gamma,\varepsilon) : \gamma \in G\}$.

\medskip

\noindent {\bf Case 2.} $\kappa \in \bigcup\{A(\alpha, \gamma) : \alpha \in I \ \wedge \ \gamma \in J_\alpha\}$.

\medskip

Let $I^{*} := \{\alpha\in I : \exists \gamma \in J_\alpha(\kappa \in A(\alpha,\gamma))\}$. When $\abs{I\setminus I^{*}} = \cf(\mu)$, then $\kappa \not\in \bigcup\{A(\alpha,\gamma) : \alpha\in I\setminus I^{*} \ \wedge \ \gamma\in J_\alpha\}$. Therefore, a similar reasoning to that presented in Case 1 produces $G\in [\mu]^{\mu}$ and $f\in C_p(L(\lambda,\kappa))$ with $f\in \bigcap\{V(g_\gamma,B_\gamma,\varepsilon) : \gamma\in G\}$.

Now suppose that $\abs{I\setminus I^{*}} < \cf(\mu)$. For each $\alpha \in I^{*}$ let $I_\alpha^{*} := \{\gamma \in J_\alpha : \kappa \in A(\alpha,\gamma)\}$. Hence, if $I' := \{\alpha \in I^{*} : \abs{I_\alpha^{*}} = \mu_{\alpha+1}\}$ satisfies $\abs{I' }<\cf(\mu)$, then $\abs{I^* \setminus I'} = \cf(\mu)$ and $\kappa \not\in \bigcup\{A(\alpha,\gamma ) : \alpha\in I^{*}\setminus I'\ \wedge \ \gamma\in J_\alpha \setminus I_{\alpha}^{*}\}$. So, once again, arguments similar to those present in Case 1 provide $G\in [\mu]^{\mu}$ and $f\in C_p(L(\lambda,\kappa))$ with $ f\in \bigcap\{V(g_\gamma,B_\gamma,\varepsilon) : \gamma\in G\}$.

It remains to consider when $\abs{I'} = \cf(\mu)$. Observe that, as $\kappa \in \bigcap\{A(\alpha,\gamma) : \alpha\in I' \ \wedge \ \gamma \in I_\alpha^*\}$, we can assume without losing generality that $\xi_0 = \kappa$. In these circumstances, a deduction similar to the one set out in the Claim of Lemma~\ref{lemma_Cp_L(lambda,kappa)_caliber_mu_1} allows us to prove the folowing statements.

\begin{enumerate}
\item For any $\alpha \in I'$ there exist collections $\{E(\alpha,i) : i<m\}$ and $\{x(\alpha,i) : i<m\}$ such that the following conditions are satisfied for every $i<m$:

\begin{enumerate}
\item $E(\alpha,0) \in [I_\alpha^{*}]^{\mu_{\alpha+1}}$,
\item if $i+1<m$, then $E(\alpha,i+1) \in [E(\alpha,i)]^{\mu_{\alpha+1}}$, and
\item $x(\alpha,i) \in \bigcap \{B_{\varepsilon/2}(f(\alpha,\gamma)(\xi_i)) : \gamma \in E(\alpha,i)\}$.
\end{enumerate}

\item There are collections $\{K_i : i<m\}$ and $\{x_i : i<m\}$ such that the following conditions are satisfied for any $i<m$:

\begin{enumerate}
\item $K_0 \in [I']^{\cf(\mu)}$,
\item if $i+1<m$, then $K_{i+1} \in [K_i]^{\cf(\mu)}$, and
\item $x_i \in \bigcap \{B_{\varepsilon/2}(x(\alpha,i)) : \alpha \in K_i\}$.
\end{enumerate}

\item For any $\alpha \in K_{m-1}$ there exist collections $\{F(\alpha,j) : j<n\}$, $\{y(\alpha,j) : j<n \}$, $\{G(\alpha,k) : k<p\}$ and $\{z(\alpha,k) : k<p\}$ such that the following conditions are satisfied for every $i<m$, $j<n$ and $k<p$:

\begin{enumerate}
\item $F(\alpha,0) \in [E(\alpha,m-1)]^{\mu _{\alpha+1}}$ and $G(\alpha,0) \in [F(\alpha,n-1)]^{\mu_{\alpha+1}}$,
\item if $j+1<n$ and $k+1<p$, then $F(\alpha,j+1) \in [F(\alpha,j)]^{\mu_{\alpha+1 }}$ and $G(\alpha,k+1) \in [G(\alpha,k)]^{\mu_{\alpha+1}}$,
\item $y(\alpha,j) \in \bigcap \{B_\varepsilon(f(\alpha,\gamma)(\xi(\alpha,j))) : \gamma \in F(\alpha,j )\}$, and
\item $z(\alpha,k) \in \bigcap \{B_\varepsilon(f(\alpha,\gamma)(\xi(\alpha,\gamma,k))) : \gamma \in G(\alpha,k)\}$.

\end{enumerate}

\end{enumerate}

Finally, consider the function $f: L(\lambda,\kappa) \to \mathbb{R}$ determined by $$f(\xi) := \begin{cases} x_0, & \text{if} \ \xi \in L(\lambda,\kappa) \setminus \bigcup\{A(\alpha,\gamma) : \alpha \in K_{m-1} \ \wedge \ \gamma\in G(\alpha,p-1)\}, \\
x_i, & \text{if} \ \exists i<m\left(\xi = \xi_i\right), \\
y(\alpha,j), & \text{if} \ \exists\alpha\in K_{m-1} \exists j<n\left(\xi = \xi(\alpha,j)\right), \\
z(\alpha,k), & \text{if} \ \exists\alpha\in K_{m-1} \exists \gamma \in G(\alpha,p-1) \exists k<p\left(\xi = \xi(\alpha,\gamma,k)\right).
\end{cases}$$ Therefore, since $\abs{\bigcup\{A(\alpha,\gamma)\setminus\{\kappa\} : \alpha \in K_{m-1} \ \wedge \ \gamma\in G( \alpha,p-1)\}} \leq \mu < \lambda$, $\kappa \not\in \bigcup\{A(\alpha,\gamma)\setminus\{\kappa\} : \alpha \in K_{m-1} \ \wedge \ \gamma\in G(\alpha,p-1)\}$ and $f[L(\lambda,\kappa) \setminus \bigcup\{A(\alpha, \gamma)\setminus\{\kappa\} : \alpha \in K_{m-1} \ \wedge \ \gamma\in G(\alpha,p-1)\}] \subseteq \{x_0\}$, $f$ is a continuous function. Lastly, notice that $G := \bigcup\{G(\alpha, p-1) : \alpha \in K_{m-1}\}$ satisfies the relation $\abs{G} = \mu$ and, furthermore, $f$ is an element of $\bigcap\{V(g_\gamma,B_\gamma,\varepsilon) : \gamma \in G\}$.
\end{proof} 

With these previous results we can precisely determine the calibers of $C_p\left(L(\lambda,\kappa)\right)$.

\begin{theorem}\label{thm_Cp_calibers_lambda_Lindelof} If $\kappa$ and $\lambda$ are cardinals with $\omega \leq \lambda\leq \kappa$, then \begin{align*} \C\left(C_p\left(L(\lambda, \kappa)\right)\right) = &\{\mu\in\CN : \cf(\mu)>\kappa\} \cup \{\mu\in\CN : \cf(\mu)>\omega \ \wedge \ \mu < \lambda\} \cup \\
&\{\mu\in\CN : \cf(\mu)>\omega \ \wedge \ \cf(\mu)<\lambda \ \wedge \ \kappa<\mu\}.
\end{align*} Equivalently, $$\C\left(C_p\left(L(\lambda, \kappa)\right)\right)=\left\{\mu\in\CN : \cf(\mu)>\omega \ \wedge \ \mu\not\in [\lambda ,\kappa] \ \wedge \ \cf(\mu)\not\in [\lambda,\kappa]\right\}.$$

\end{theorem}

\begin{proof} Let $A := \{\mu\in\CN : \cf(\mu)>\kappa\}$, $B := \{\mu\in\CN : \cf(\mu)>\omega \ \wedge \ \mu < \lambda\}$ and $C := \{\mu\in\CN : \cf(\mu)>\omega \ \wedge \ \cf(\mu)<\lambda \ \wedge \ \kappa<\mu\}$. To check that all the elements of $A\cup B \cup C$ are calibers for $C_p\left(L(\lambda,\kappa)\right)$, note that $L(\lambda,\kappa)$ condenses onto $L(\omega,\kappa)$. For this reason, Propositions~\ref{prop_Cp_condensation} and \ref{prop_Cp_iw_Tychonoff} combined with the equality $nw\left(L(\omega,\kappa)\right)=\kappa$ imply that $$A\subseteq \ C\left(C_p\left(L(\omega,\kappa)\right)\right)\subseteq \C\left(C_p\left(L(\lambda,\kappa)\right)\right).$$ Furthermore, Lemma~\ref{lemma_Cp_L(lambda,kappa)_caliber_mu_1}, Corollary~\ref{cor_Cp_L(lambda,kappa)_pi_peso} and Lemma~\ref{lemma_Cp_L(lambda,kappa)_caliber_mu_2} guarantee that the elements of $B\cup C$ are calibers for $C_p\left(L(\lambda, \kappa)\right)$.

For the remaining inclusion, assume that $\mu\in \CN$ is not an element of $A\cup B \cup C$. Since the equality $\cf(\mu) = \omega$ implies that $\mu$ is not a caliber for $C_p\left(L(\lambda, \kappa)\right)$ by Proposition~\ref{cor_cellularity_cofinality}, we can additionally assume that $\cf(\mu)>\omega$. Furthermore, notice that the relations $\cf(\mu) \leq \kappa$ and $\lambda \leq \mu$ hold because $\mu\not \in A \cup B$.

It remains to analyze what happens when $\kappa<\mu$ and when $\mu \leq \kappa$. If $\kappa<\mu$, the condition $\mu\not \in C$ ensures that $\lambda \leq \cf(\mu)$ and, therefore, the relation $\cf(\mu) \leq \kappa$ and Corollary~\ref{cor_Cp_L(lambda,kappa)_no_caliber_mu} guarantee that $\mu$ is not a caliber for $C_p\left(L(\lambda, \kappa)\right)$. Finally, if $\mu \leq \kappa$, the relation $\lambda \leq \mu$ and Lemma~\ref{lemma_Cp_L(lambda,kappa)_caliber_mu_1} confirm that $\mu$ is not a caliber for $C_p \left(L(\lambda, \kappa)\right)$.
\end{proof}

As a last remark, we mention that Comfort and Negrepontis obtained the following result related to the calibers of $C_p\left(L(\lambda,\kappa),D(2)\right)$ (see \cite[Theorem~4.2, p.~80]{comneg1982}):

\begin{theorem} If $\kappa$ and $\lambda$ are cardinals with $\omega\leq \lambda \leq \kappa$, then $$\C\left(C_p\left(L(\lambda, \kappa),D(2)\right)\right)=\left\{\mu\in\CN : \cf(\mu)>\omega \ \wedge \ \mu\not\in [\lambda ,\kappa] \ \wedge \ \cf(\mu)\not\in [\lambda,\kappa]\right\}.$$

\end{theorem}

%%%%%%%%%%%%%%%%%%%%%%%%%%%%%%%%%%%%%%%%%%%%%%%%%%%%%%%%%%%%%%
\section{Calibers for $C_p(X)$ when $X$ is a topological product}\label{secc_productos}
%%%%%%%%%%%%%%%%%%%%%%%%%%%%%%%%%%%%%%%%%%%%%%%%%%%%%%%%%%%%%%

Regarding the calibers of $C_p$-spaces over topological products, Proposition~\ref{prop_P(kappa,2)_prod} and Corollary~\ref{cor_P(kappa,2)_prod} were our initial steps in this direction. To obtain more results we first have to establish a generalization of \cite[Lemma~9, p.~683]{bailey2007} to calculate the $i$-weight of an arbitrary product.

\begin{proposition}\label{prop_iw_products} If $\lambda$ is a cardinal number and $\{X_\alpha : \alpha<\lambda\}$ is a collection of non-trivial Tychonoff spaces, then $$iw\left(\prod_{\alpha<\lambda} X_\alpha\right)=\lambda\cdot\sup\{iw(X_\alpha) : \alpha<\lambda\}.$$
\end{proposition}

\begin{proof} Let $X:=\prod_{\alpha<\lambda} X_\alpha$ and $\mu:=\lambda\cdot\sup\{iw(X_\alpha) : \alpha<\lambda\}$. If for each $\alpha<\lambda$ we let $Y_\alpha$ be a Tychonoff space such that $X_\alpha$ is condensed onto $Y_\alpha$ and $iw(X_\alpha)=w(Y_\alpha)$, then $X$ condenses onto the Tychonoff space $\prod\{Y_\alpha : \alpha<\lambda\}$; consequently, $$iw\left(X\right) \leq w\left(\prod_{\alpha<\lambda} Y_\alpha\right) = \lambda\cdot\sup\{w(Y_\alpha) : \alpha<\lambda\} = \mu.$$

To determine the opposite inequality observe that if $Y$ is a Tychonoff space such that $X$ condenses onto $Y$ and $iw(X)=w(Y)$, then, as each $X_\alpha$ embeds into the product $X$, we infer that all $X_\alpha$ are condensed onto a subspace of $Y$; therefore, we obtain that $\sup\{iw(X_\alpha) : \alpha<\lambda\} \leq w(Y) = iw(X)$.

What remains is to notice that if $1\leq \lambda<\omega$, what has been done so far guarantees that $iw(X) = \mu$. If $\lambda\geq \omega$, we use that the Cantor cube $D(2)^{\lambda}$ embeds into $X$ and that the $i$-weight is monotonic to ensure that $$\lambda = nw\left(D(2)^{\lambda}\right) = iw\left(D(2)^{\lambda}\right) \leq iw\left(X\right).$$ Thus, $\mu\leq iw(X)$.
\end{proof}

A consequence of Proposition~\ref{prop_iw_products} and the Hewitt-Marczewski-Pondiczery Theorem (see \cite[Theorem~11.2, p~42]{hodel1984}) is that the countable $i$-weight is preserved under countable products and separability is preserved under products of at most $\mathfrak{c}$ factors. These facts and Corollary~\ref{cor_prop_Cp_iw_Tychonoff} imply the following.

\begin{corollary} Let $1\leq n<\omega$ and $\{X_\alpha : \alpha<\lambda\}$ be a family of Tychonoff spaces. If we denote by $X$ the topological product $\prod_{\alpha<\lambda} X_\alpha$, then the following statements are true.
\begin{enumerate}
\item If each $X_\alpha$ has countable $i$-weight and $1\leq \lambda\leq\omega$, then $C_{p,2n-1}(X)$ is strongly \v{S}anin.
\item When all $X_\alpha$ are separable and $1\leq \lambda\leq \mathfrak{c}$, $C_{p,2n}(X)$ is strongly \v{S}anin.
\item Consequently, if each $X_\alpha$ is a separable space with countable $i$-weight and $1\leq \lambda\leq\omega$, then $C_{p,n}(X)$ is strongly \v{S}anin.
\end{enumerate}

\end{corollary}

For example, if we denote the Sorgenfrey line by $\mathbb{S}$, then $C_{p,n }\left(\mathbb{S}^{\lambda}\right)$ is strongly \v{S}anin for every $1\leq \lambda\leq\omega$ and $1\leq n< \omega$.

Furthermore, with Propositions~\ref{prop_Cp_iw_Tychonoff} and~\ref{prop_iw_products} we can determine some of the calibers for the $C_p$-space of the products of non-trivial Tychonoff spaces:

\begin{proposition}\label{prop_calibers_products_1} Let $\lambda$ be a cardinal number and $\{X_\alpha : \alpha<\lambda\}$ a collection of non-trivial Tychonoff spaces. If $\mu:= \lambda\cdot\sup\{iw(X_\alpha) : \alpha<\lambda\}$ and $X_\pi := \prod_{\alpha<\lambda} X_\alpha$, then $$\{\kappa\in\CN : \cf(\kappa)>\mu\} \subseteq \C\left(C_p\left(X_\pi\right)\right).$$ In particular, when $X$ is a non-trivial Tychonoff space, $$\{\kappa\in\CN : \cf(\kappa)>\lambda\cdot iw(X)\} \subseteq \C\left(C_p\left(X^{\lambda}\right)\right).$$
\end{proposition}

Recall that if $\{X_\alpha : \alpha<\lambda\}$ is a family of topological spaces and $Y$ is also a space, a function $f:\prod_{\alpha<\lambda} X_\alpha \to Y$ {\it depends only} on $J\subseteq \lambda$ if for any $x,y\in \prod_{\alpha<\lambda} X_\alpha$, the equality $\restr{x}{J} =\restr{y}{J}$ implies that $f(x)=f(y)$. Furthermore, $f$ {\it depends on less than $\kappa$ coordinates} if there exists $J\in [\lambda]^{<\kappa}$ such that $f$ depends only on $J$.

A factorization result due to Mibu (see \cite{mibu1944}) states that if $\lambda$ is a cardinal number, $\{X_\alpha : \alpha<\lambda\}$ is a familiy of compact Hausdorff spaces and $f:\prod_{\alpha<\lambda} X_\alpha \to \mathbb{R}$ is a continuous function, then $f$ depends on less than $\omega_1$ coordinates.

Now we need to determine the collection $\C(C_p(D(2)^{\lambda}))$ for each cardinal number $\lambda\geq \omega$. A consequence of Mibu's theorem is that if $\lambda$ is a cardinal and $f\in C\left(D(2)^{\lambda}\right)$, then $f$ depends on less than $\omega_1$ coordinates. The following result is supported by this fact.

\begin{proposition}\label{prop_Cp_Cantor_no_caliber_kappa} If $\lambda$ is an infinite cardinal, then $C_p(D(2)^{\lambda})$ does not have caliber $\lambda$.
\end{proposition}

\begin{proof} First, if $\cf(\lambda)=\omega$, then Proposition~\ref{cor_cellularity_cofinality} ensures that $\lambda$ is not a caliber for $C_p\left(D(2)^{\lambda}\right )$. When $\cf(\lambda)>\omega$ we do the following: for each $\alpha<\lambda$, let $x_\alpha : \lambda \to D(2)$ be the characteristic function of $\{ \alpha\}$; let $U_\alpha$ be a clopen subset of $D(2)^{\lambda}$ with the property that $\{x_\alpha\} = U_\alpha \cap \{x_\beta: \beta<\lambda\}$, and we denote by $\chi_\alpha : D(2)^{\lambda} \to \mathbb{R}$ the characteristic function of $U_\alpha$.

Consider the collection $\left\{V(\chi_\alpha, \{x_{\alpha},x_{\alpha+1}\}, 1/2) : \alpha<\lambda\right\}$ and take $J\in [\lambda]^{\lambda}$. In search of a contradiction, suppose there exists a function $$f\in \bigcap\left\{V(\chi_\alpha, \{x_{\alpha},x_{\alpha+1}\}, 1/2 ) : \alpha\in J\right\}.$$ Let $A\in [\lambda]^{\leq\omega}$ be such that, for any $x,y\in D(2)^{\lambda }$, $f(x)=f(y)$ whenever $\restr{x}{A}=\restr{y}{A}$. Finally, if $\alpha \in J \cap (\sup A, \lambda)$, then $\restr{x_\alpha}{A}=\restr{x_{\alpha+1}}{A}$ and, therefore, if we set $r:=f(x_\alpha)=f(x_{\alpha+1})$, the condition $f\in V(\chi_\alpha, \{x_{\alpha}, x_{\alpha+1}\}, 1/2)$ implies the absurdity $r<1/2$ and $1/2<r$.
\end{proof}

\begin{proposition}\label{prop_Cp_calibers_Cantor_cube} If $\lambda$ is an infinite cardinal, then $$\C\left(C_p\left(D(2)^{\lambda}\right)\right)=\{\kappa\in\CN : \cf(\kappa)>\lambda\}.$$
\end{proposition}

\begin{proof} To begin with, since $nw\left(D(2)^{\lambda}\right)=\lambda$, Proposition~\ref{prop_Cp_iw_Tychonoff} ensures that $$\{\kappa\in\CN : \cf( \kappa)>\lambda\}\subseteq \C\left(C_p\left(D(2)^{\lambda}\right)\right).$$ Besides, if $\kappa\in\CN$ satisfies $\cf(\kappa)\leq\lambda$, then a combination of Propositions~\ref{prop_Cp_subspace} and~\ref{prop_Cp_Cantor_no_caliber_kappa} implies that $C_p(D(2)^{\lambda})$ does not have caliber $\cf(\kappa)$. Finally, Proposition~\ref{prop_caliber_cofinality} guarantees that $\kappa$ is not a caliber for $C_p(D(2)^{\lambda})$.
\end{proof}

From the previous result we can precisely determine the collection of calibers of the $C_p$ of a topological product, provided that each of the factor spaces satisfies two simple inequalities.

\begin{proposition}\label{prop_Cp_calibers_product} Let $\lambda$ be an infinite cardinal and $\{X_\alpha : \alpha<\lambda\}$ a family of non-trivial Tychonoff spaces such that $iw(X_\alpha)\leq \lambda$ for every $\alpha<\lambda$. If $X_\pi := \prod_{\alpha<\lambda} X_\alpha$, then $$\C\left(C_p\left(X_\pi\right)\right)=\{\kappa\in\CN : \cf(\kappa)>\lambda\}.$$ 
Hence, if $X$ is a non-trivial Tychonoff space with $iw(X)\leq \lambda$, then $\C(C_p(X^{\lambda}))=\{\kappa\in\CN : \cf(\kappa)>\lambda\}$.
\end{proposition}

\begin{proof} Propositions~\ref{prop_Cp_subspace}, \ref{prop_calibers_products_1} and \ref{prop_Cp_calibers_Cantor_cube} ensure, respectively, that \begin{align*}
\C\left(C_p\left(X_\pi\right)\right) &\subseteq \C\left(C_p\left(D(2)^{\lambda}\right)\right),  \\
\{\kappa\in\CN : \cf(\kappa)>\lambda\} &\subseteq \C\left(C_p\left(X_\pi\right)\right) \ \text{and} \\
 \C\left(C_p\left(D(2)^{\lambda}\right)\right) &\subseteq \{\kappa\in\CN : \cf(\kappa)>\lambda\}.
\end{align*} Thus, $\C\left(C_p(X_\pi)\right)=\{\kappa\in\CN : \cf(\kappa)>\lambda\}$.
\end{proof}

Arhangel'ski\u{\i} and Tkachuk showed in \cite{arhtka1986} that if $X$ is a compact Hausdorff space and $\cf\left(w(X)\right)>\omega$, then $ \cf\left(w(X)\right)$ is not a caliber for $C_p(X)$. This result is useful to verify that the restriction on the $i$-weight of the space is necessary in Proposition~\ref{prop_Cp_calibers_product}:

\begin{proposition} If $X$ is a non-trivial compact Hausdorff space and $\lambda$ is an infinite cardinal with $\cf(w(X))>\lambda$, then $$\{w(X),\cf(w(X))\} \subseteq \{\kappa\in\CN : \cf(\kappa)>\lambda\}\setminus \C\left(C_p\left(X^{\lambda}\right)\right).$$

\end{proposition}

\begin{proof} Clearly, both $w(X)$ and $\cf(w(X))$ belong to $\{\kappa\in\CN : \cf(\kappa)>\lambda\}$. On the other hand, $w(X^{\lambda})=\lambda\cdot w(X) = w(X)$ (see \cite[Ch.~5]{juhasz1980}); in particular, $X^{\lambda}$ is a compact Hausdorff space with $\cf(w(X^{\lambda}))=\cf(w(X))>\lambda\geq\omega$ and, therefore, $\cf(w(X))$ is not a caliber for $C_p(X^\lambda)$ by the theorem of Arhangel'ski\u{\i} and Tkachuk. Finally, Proposition~\ref{prop_caliber_cofinality} guarantees that $w(X)$ is also not a caliber for $C_p(X^\lambda)$.
\end{proof}

For example, for every infinite cardinal $\kappa$, $\kappa^{+}$ is an element of $\{\mu\in \CN : \cf(\mu)>\kappa\}$ and is not a caliber for $C_p([0,\kappa^+]^{\kappa})$. Thus, the restriction {\lq\lq}$iw(X)\leq \lambda${\rq\rq} cannot be dropped in Proposition~\ref{prop_Cp_calibers_product}.

%%%%%%%%%%%%%%%%%%%%%%%%%%%%%%%%%%%%%%%%%%%%%%%%%%%%%%%%%%%%%%
\section{Calibers for $C_p(X)$ when $X$ is a topological sum}\label{secc_sumas}
%%%%%%%%%%%%%%%%%%%%%%%%%%%%%%%%%%%%%%%%%%%%%%%%%%%%%%%%%%%%%%

For a family of topological spaces $\{X_\alpha : \alpha<\lambda\}$ we use the symbol $X_s$ to denote the topological sum $\bigoplus_{\alpha<\lambda} X_\alpha$. The purpose of this section is to obtain results regarding the set $\C(C_p(X_s))$.

A first remark is that, by a classical duality theorem (see \cite{arh1992}), $C_p(X_s)$ is homemorphic to the topological product $\prod_{\alpha<\lambda} C_p(X_\alpha)$. Therefore, the known results with reference to the preservation of calibers in products are useful in determining the chain conditions of $C_p(X_s)$. The following corollary is a consequence of Theorem~\ref{thm_RSS}.

\begin{corollary}\label{cor_Cp_sums} Let $\kappa$ be an infinite cardinal, $\lambda$ a cardinal number, $\{X_\alpha : \alpha<\lambda\}$ a family of topological spaces and $X$ a topological space.

\begin{enumerate}
\item $\UR\cap \bigcap_{\alpha<\lambda} \C(C_p(X_\alpha))\subseteq \C(C_p(X_s)) \subseteq \bigcap_{\alpha<\lambda} \C(C_p(X_\alpha))$.
\item If $\lambda<\cf(\kappa)$ and $\kappa \in \bigcap_{\alpha<\lambda} \C(C_p(X_\alpha))$, then $\kappa\in \C(C_p(X_s))$.
\item $\C(C_p(X)) = \C(C_p(X\times D(\lambda)))$.

\end{enumerate}

\end{corollary}

A natural question arises. Under what conditions is the following satisfied? \begin{align}\label{eq_C(Cp)} \C(C_p(X_s)) = \bigcap_{\alpha<\lambda} \C(C_p(X_\alpha)).
\end{align} For certain classes of spaces it is not difficult to verify that the relation (\ref{eq_C(Cp)}) always holds.

\begin{proposition}\label{prop_Cp_sum_discrete} If $\lambda$ is a cardinal number and $\{X_\alpha : \alpha<\lambda\}$ is a family of discrete spaces, then $$\C(C_p(X_s)) = \bigcap_{\alpha<\lambda} \C(C_p(X_\alpha)).$$

\end{proposition}

\begin{proof} If $\mu := \lambda \cdot \sup\left\{\abs{X_\alpha} : \alpha<\lambda\right\}$, then $$C_p(X_s) \cong \prod_{\alpha<\lambda} C_p(X_\alpha) \cong \prod_{\alpha<\lambda} \mathbb{R}^{X_\alpha} \cong \mathbb{R}^{\mu}$$ and, therefore, Proposition~\ref{prop_Cp_precalibers} implies that $$\C(C_p(X_s)) = \C\left(\mathbb{R}^{\mu}\right) = \UC = \bigcap_{\alpha<\lambda} \C\left(\mathbb{R}^{X_\alpha}\right)=\bigcap_{\alpha<\lambda} \C(C_p(X_\alpha)).$$
\end{proof}

Now, to provide a result in this sense for the particular spaces that we presented in the previous sections, what follows is to prove a new consistency theorem regarding the preservation of chain conditions in topological products. The first thing we do is gather some auxiliary results. The following equivalence is well known.

\begin{lemma}\label{lemma_weakly_inaccesible} The following statements are equivalent for a regular and uncountable cardinal $\kappa$.

\begin{enumerate}
\item $\kappa$ is a limit cardinal (equivalently, $\kappa$ is weakly inaccessible).
\item $\abs{\{\lambda \in \CN : \lambda<\kappa\}} = \kappa$.
\end{enumerate}

\end{lemma}

\begin{lemma}\label{lemma_homeomorphic_product} Let $\kappa$ be an infinite cardinal with $\cf(\kappa)>\omega$, $\lambda$ a cardinal number and $\{X_\alpha : \alpha<\lambda\}$ a family of topological spaces. Let $\mu<\cf(\kappa)$ be a cardinal number and $\{Y_\beta : \beta<\mu\}$ be a family of topological spaces that satisfy the following conditions:

\begin{enumerate}
\item for any $\alpha<\lambda$, there exists $\beta<\mu$ with $X_\alpha \cong Y_\beta$,
\item for all $\beta<\mu$, there is $\alpha<\lambda$ such that $Y_\beta \cong X_\alpha$, and
\item when $\beta<\gamma<\mu$, $Y_\beta$ is not homeomorphic to $Y_\gamma$.
\end{enumerate}

Then, if $\A\in \left\{\C, \P, \WP\right\}$, $\kappa \in \bigcap_{\alpha<\lambda} \A (X_\alpha)$ implies $\kappa \in \A(\prod_{\alpha<\lambda} X_\alpha)$.

\end{lemma}

\begin{proof} If for each $\beta<\mu$ we define $J_\beta := \{\alpha<\lambda : X_\alpha \cong Y_\beta\}$ and $Z_\beta := Y_\beta^{|J_ \beta|}$, then clearly $\{J_\beta : \beta<\mu\}$ is a partition of $\lambda$ and $$\prod_{\alpha<\lambda} X_\alpha \cong \prod_{\beta<\mu} \left(\prod_{\alpha\in J_\beta} X_\alpha\right) \cong \prod_{\beta<\mu} Z_\beta.$$ Thus, since $\kappa\in \bigcap_{\beta<\mu} \A(Y_\beta)$, by Theorem~\ref{thm_RSS}(3) it is satisfied that $\kappa\in \A(Z_\beta)$ whenever $\beta<\mu$ and, therefore, Theorem~\ref{thm_RSS}(2) implies that $\kappa \in \A\left(\prod_ {\beta<\mu} Z_\beta\right)$.
\end{proof}

\begin{proposition}\label{prop_equivalence_calibers} The following statements are equivalent for $\A\in \left\{\C, \P, \WP\right\}$.

\begin{enumerate}

\item If $\lambda$ is a cardinal number, $\{\lambda_\alpha : \alpha<\lambda\}$ is a subset of $\mathsf{CN}$ and $\{X_\alpha : \alpha< \lambda\}$ is a family of topological spaces that satisfy the following conditions for any $\alpha<\lambda$ and $\beta<\lambda$:

\begin{enumerate}
\item $\A(X_\alpha) = \{\kappa \in \mathsf{CN} : \cf(\kappa)>\lambda_\alpha\}$, and
\item if $\lambda_\alpha = \lambda_\beta$, then $X_\alpha \cong X_\beta$,
\end{enumerate} then $\bigcap_{\alpha<\lambda} \A(X_\alpha) = \A(\prod_{\alpha<\lambda} X_\alpha)$.

\item If $\lambda$ is a cardinal number, $\{\lambda_\alpha : \alpha<\lambda\}$ is a subset of $\mathsf{CN}$ and $\{X_\alpha : \alpha< \lambda\}$ is a family of topological spaces that satisfy the following conditions for any $\alpha<\lambda$:

\begin{enumerate}
\item $\A(X_\alpha) = \{\kappa \in \mathsf{CN} : \cf(\kappa)>\lambda_\alpha\}$, and
\item $\{\lambda_\beta : \beta<\lambda\}$ is strictly increasing,
\end{enumerate} then $\bigcap_{\alpha<\lambda} \A(X_\alpha) = \A(\prod_{\alpha<\lambda} X_\alpha)$.

\end{enumerate}

\end{proposition}

\begin{proof} To see that (1) implies (2), suppose that conditions (2)(a) and (2)(b) are true. We need to show that (1)(b) is true as well. Indeed, if $\alpha<\lambda$ and $\beta<\lambda$ satisfy $\lambda_\alpha = \lambda_\beta$, then (2)(b) implies that $\alpha=\beta$ and, therefore, $X_\alpha = X_\beta$; in particular, $X_\alpha \cong X_\beta$.

For the reverse implication, suppose that conditions (1)(a) and (1)(b) are true. Let $\{\mu_\beta : \beta<\mu\}$ be a strictly increasing enumeration of $\{\lambda_\alpha : \alpha<\lambda\}$. For each $\beta<\mu$, let $J_\beta := \{\alpha < \lambda : \lambda_\alpha = \mu_\beta\}$. Furthermore, let $f : \mu \to \lambda$ be a function such that $f(\beta) \in J_\beta$ for all $\beta<\mu$. Clearly, $\{J_\beta : \beta <\mu\}$ is a partition of $\lambda$.

Now, if $\beta<\mu$, then $$\A(X_{f(\beta)}) = \left\{\kappa \in \mathsf{CN} : \cf(\kappa)>\lambda _{f(\beta)}\right\}=\left\{\kappa \in \mathsf{CN} : \cf(\kappa)>\mu_\beta\right\} .$$ For each $\beta<\mu$, set $Z_\beta := (X_{f(\beta)})^{|J_\beta|}$. It then turns out that $$\A(Z_\beta) =\left\{\kappa \in \mathsf{CN} : \cf(\kappa)>\mu_\beta\right\}$$ by Theorem~\ref{thm_RSS}(3). Consequently, (2) implies that $\bigcap_{\beta<\mu} \A(Z_\beta) = \A(\prod_{\beta<\mu} Z_\beta)$. Finally, since $$\prod_{\alpha<\lambda} X_\alpha \cong \prod_{\beta<\mu} \left(\prod_{\alpha\in J_\beta} X_\alpha\right) \cong \prod_{\beta<\mu} Z_\beta \quad \text{and} \quad \bigcap_{\alpha<\lambda} \A(X_\alpha) = \bigcap_{\beta<\mu} \A(Z_\beta),$$ we obtain $\bigcap_{\alpha<\lambda} \A(X_\alpha) = \A(\prod_{\alpha<\lambda} X_\alpha)$.
\end{proof}

It turns out that for families of spaces like those mentioned in Proposition~\ref{prop_equivalence_calibers}, the absence of cardinals in $\bigcap_{\alpha<\lambda} \A(X_\alpha)$ with weakly inaccessible cofinality is enough to guarantee that $\bigcap_{\alpha<\lambda} \A(X_\alpha) = \A(\prod_{\alpha<\lambda} X_\alpha)$.

\begin{theorem}\label{thm_weakly_in_calibers} Let $\A\in \left\{\C, \P, \WP\right\}$. Furthermore, let $\{\lambda_\alpha : \alpha<\lambda\}$ and $\{X_\alpha : \alpha< \lambda\}$ be as in Proposition~\ref{prop_equivalence_calibers}(1). If $\cf(\kappa)$ is not weakly inaccessible for every $\kappa \in \bigcap_{\alpha<\lambda} \A(X_\alpha)$, then $\bigcap_{\alpha<\lambda} \A(X_\alpha) = \A(\prod_{\alpha<\lambda} X_\alpha)$.

\end{theorem}

\begin{proof} Observe that, since $\bigcap_{\alpha<\lambda} \A(X_\alpha) \supseteq \A(\prod_{\alpha<\lambda} X_\alpha)$, we only need to prove $\bigcap_{\alpha<\lambda} \A(X_\alpha) \subseteq \A(\prod_{\alpha<\lambda} X_\alpha)$. Let $\kappa$ be an element of $\bigcap_{\alpha<\lambda} \A(X_\alpha)$ and suppose that $\kappa$ does not belong to $\A(\prod_{\alpha<\lambda} X_\alpha)$.

Let $\sim$ be the equivalence relation in $\lambda$ defined by $\alpha \sim \beta$ if and only if $X_\alpha \cong X_\beta$. If $J$ is a complete system of representatives for $\sim$, $\{X_\beta : \beta \in J\}$ satisfies conditions (1), (2) and (3) of Lemma~\ref{lemma_homeomorphic_product}. Therefore, the relation $\kappa \in \bigcap_{\alpha<\lambda} \A(X_\alpha) \setminus \A(\prod_{\alpha<\lambda} X_\alpha)$ implies that $\abs{J} \geq \cf(\kappa)$.

Now notice two things: first, the condition $\kappa \in \bigcap_{\alpha<\lambda} \A(X_\alpha)$ guarantees that $\{\lambda_\alpha : \alpha<\lambda\} \subseteq \cf(\kappa)$; second, if $\alpha,\beta \in J$ satisfy $\alpha\neq \beta$, then $\lambda_\alpha \neq \lambda_\beta$ because $X_\alpha$ is not homeomorphic to $X_\beta $. Hence, $\{\lambda_\alpha : \alpha\in J\}$ is a subset of $\cf(\kappa)$ of cardinality $\cf(\kappa)$ and, therefore, $\cf (\kappa)$ is weakly inaccessible by Lemma~\ref{lemma_weakly_inaccesible}.
\end{proof}

\begin{corollary} If there are no weakly inaccesible cardinals, then both statements in Proposition~\ref{prop_equivalence_calibers} hold.

\end{corollary}

\begin{corollary}\label{cor_Cp_inaccesible} Let $\lambda$ be a cardinal number, $\{X_\alpha : \alpha<\lambda\}$ a family of spaces and $\{\lambda_\alpha : \alpha<\lambda\}$ a collection of cardinals such that the following conditions hold for any $\alpha<\lambda$ and $\beta<\lambda$:

\begin{enumerate}
\item $\C(C_p(X_\alpha)) = \{\kappa \in \CN : \cf(\kappa) > \lambda_\alpha\}$, and
\item if $\lambda_\alpha = \lambda_\beta$, then $C_p(X_\alpha) \cong C_p(X_\beta)$.
\end{enumerate} If $\cf(\kappa)$ is not weakly inaccessible for every $\kappa \in \bigcap_{\alpha<\lambda} \C(C_p(X_\alpha))$, then $$\C(C_p(X_s)) = \bigcap_{\alpha<\lambda} \C(C_p(X_\alpha)).$$

\end{corollary}
 
It is important to highlight that by Proposition~\ref{prop_Cp_calibers_Cantor_cube}, if each $X_\alpha$ is of the form $D(2)^{\kappa_\alpha}$, then conditions (1) and (2) established in Corollary~\ref{cor_Cp_inaccesible} are satisfied. Hence, if there are no weakly inaccessible cardinals, the relation~(\ref{eq_C(Cp)}) holds.

However, for certain particular cases the hypothesis about the non-existence of weakly inaccessible cardinals is not necessary:

\begin{lemma}\label{lemma_Cp_cardinality} Let $\lambda$ be a cardinal number and $\{X_\alpha : \alpha<\lambda\}$ be a family of spaces such that the following conditions hold for every $\alpha<\lambda$ and $\beta <\lambda$:

\begin{enumerate}
\item $\C(C_p(X_\alpha)) = \{\kappa \in \CN : \cf(\kappa) > |X_\alpha|\}$ and
\item $\alpha<\beta<\lambda$ implies $|X_\alpha| \neq |X_\beta|$.

\end{enumerate} Under these hypotheses, $$\C(C_p(X_s)) = \bigcap_{\alpha<\lambda} \C(C_p(X_\alpha)).$$

\end{lemma}

\begin{proof} Let $\kappa_\alpha := |X_\alpha|$ for each $\alpha<\lambda$. Note that if $\kappa \in \bigcap_{\alpha<\lambda} \C(C_p(X_\alpha))$, then $\cf(\kappa) \in \bigcap_{\alpha<\lambda} \C (C_p(X_\alpha))$ (see Proposition~\ref{prop_caliber_cofinality}) and, therefore, $\cf(\kappa) \in \C(C_p(X_s))$ by Corollary~\ref{cor_Cp_sums}(1). Furthermore, since $\{\kappa_\alpha : \alpha<\lambda\} \subseteq \cf(\kappa)$, it turns out that $\lambda = |\{\kappa_\alpha : \alpha<\lambda\}| \leq \cf(\kappa)$ and $\sup \{\kappa_\alpha : \alpha<\lambda\} \leq \cf(\kappa)$. Finally, the relations $$\pi w\left(C_p(X_s)\right) = \abs{X_s} = \lambda \cdot \sup\left\{\kappa_\alpha : \alpha<\lambda\right\} \leq \cf(\kappa) \leq \kappa,$$ imply that $\kappa \in \C(C_p(X_s))$ by virtue of Proposition~\ref{prop_caliber_pi_weight}.
\end{proof}

\begin{proposition}\label{prop_Cp_cardinality} Let $\lambda$ be a cardinal number and $\{X_\alpha : \alpha<\lambda\}$ be a family of spaces such that the following conditions hold for every $\alpha<\lambda$ and $\beta <\lambda$:

\begin{enumerate}
\item $\C(C_p(X_\alpha)) = \{\kappa \in \CN : \cf(\kappa) > |X_\alpha|\}$;
\item $X_\alpha \cong X_\beta$ whenever $|X_\alpha| = |X_\beta|$, and
\item if $I_\alpha := \{\gamma<\lambda : |X_\gamma| = |X_\alpha|\}$, then $|X_\alpha|\cdot |I_\alpha| = |X_\alpha|$.

\end{enumerate} Under these hypotheses, $$\C(C_p(X_s)) = \bigcap_{\alpha<\lambda} \C(C_p(X_\alpha)).$$

\end{proposition}

\begin{proof} Let $\lambda_\alpha := |I_\alpha|$ and $\kappa_\alpha := |X_\alpha|$ for each $\alpha<\lambda$. Let $\sim$ be the equivalence relation in $\lambda$ determined by the rule $\alpha\sim \beta$ if and only if $\kappa_\alpha = \kappa_\beta$. Note that if $J$ is a complete system of representatives for $\sim$ and $X_J := \bigoplus _{\alpha \in J} (X_\alpha \times D(\lambda_\alpha))$, then $C_p (X_s) \cong C_p(X_J)$. Furthermore, for every $\alpha \in J$ Corollary~\ref{cor_Cp_sums}(3) and item (3) ensure that \begin{align*} \C(C_p(X_\alpha \times D(\lambda_ \alpha))) &= \C(C_p(X_\alpha)) = \{\kappa \in \CN : \cf(\kappa) > \kappa_\alpha\} \\
&= \{\kappa \in \CN : \cf(\kappa) > \kappa_\alpha \cdot \lambda_\alpha\} \\
&= \{\kappa \in \CN : \cf(\kappa) > |X_\alpha \times D(\lambda_\alpha)|\}.
\end{align*} Also, if $\alpha,\beta\in J$ are distinct, then $|X_\alpha \times D(\lambda_\alpha)| = \kappa_\alpha \cdot \lambda_\alpha = \kappa_\alpha \neq \kappa_\beta = \kappa_\beta\cdot \lambda_\beta = |X_\beta\times D(\lambda_\beta)|$. Therefore, Corollary~\ref{cor_Cp_sums} and Lemma~\ref{lemma_Cp_cardinality} guarantee that \begin{align*} \bigcap_{\alpha<\lambda} \C(C_p(X_\alpha)) &\subseteq \bigcap _{\alpha\in J} \C(C_p(X_\alpha)) = \bigcap _{\alpha\in J} \C(C_p(X_\alpha \times D(\lambda_\alpha))) \\
&= \C(C_p(X_J)) = \C(C_p(X_s)) \subseteq \bigcap_{\alpha<\lambda} \C(C_p(X_\alpha).
\end{align*}
\end{proof}

The following corollary is a consequence of Proposition~\ref{prop_calibers_Cp_[0,kappa]}, Theorem~\ref{thm_Cp_calibers_lambda_Lindelof} and Proposition~\ref{prop_Cp_cardinality}.

\begin{corollary} Let $\lambda$ be a cardinal number and $\{\kappa_\alpha : \alpha<\lambda\}$ be a family of distinct cardinals. If each $X_\alpha$ is of the form $[0,\kappa_\alpha]$, $[0,\kappa_\alpha)$, $L(\omega,\kappa_\alpha)$ or $L(\omega_1,\kappa_\alpha)$, then $$\C(C_p(X_s)) = \bigcap_{\alpha<\lambda} \C(C_p(X_\alpha)).$$

\end{corollary}

To conclude this work we present an example to verify that, in general, equality~(\ref{eq_C(Cp)}) is not necessarily true. We recommend that the reader consults Section~\ref{secc_Lindelof} to recall the notation defined there.

\begin{proposition} Let $\kappa$ and $\lambda$ be cardinals with $\kappa>\lambda>\cf(\kappa)>\omega$ and $\cf(\lambda)>\omega$. Furthermore, let $\{\kappa_\alpha : \alpha<\cf(\kappa)\}$ be a family of cardinals with the following characteristics:
\begin{enumerate}
\item $\lambda \leq \kappa_0$,
\item $\{\kappa_\alpha : \alpha<\cf(\kappa)\}$ is strictly increasing, and
\item $\sup\{\kappa_\alpha : \alpha<\cf(\kappa)\}=\kappa$.
\end{enumerate} With this background, $$\kappa \in \bigcap_{\alpha<\cf(\kappa)} \C\left(C_p\left(L(\lambda, \kappa_\alpha)\right) \right) \setminus \C\left(C_p\left(\bigoplus_{\alpha<\cf(\kappa)} L(\lambda, \kappa_\alpha)\right)\right).$$

\end{proposition}

\begin{proof} For each $\alpha<\cf(\kappa)$ let $X_\alpha :=L(\lambda, \kappa_\alpha)\times \{\alpha\}$. Since $\cf(\kappa)>\omega$, $\cf(\kappa) < \lambda$ and $\kappa_\alpha<\kappa$ for any $\alpha<\cf(\kappa)$, $\kappa$ is a caliber for each $C_p(X_\alpha)$ by Theorem~\ref{thm_Cp_calibers_lambda_Lindelof}. To verify that $\kappa$ is not a caliber for $C_p(X_s)$, let $J_0 := [0,\kappa_0)$ and $J_\alpha := [\kappa_\alpha, \kappa_{\alpha+ 1})$ whenever $\alpha<\cf(\kappa)$. Furthermore, for all $\xi < \kappa$ let $\alpha(\xi)<\cf(\kappa)$ be the sole ordinal such that $\xi \in J_{\alpha(\xi)}$.

Now, for each $\xi<\kappa$ let $f_\xi : X_s \to \mathbb{R}$ be the continuous function given by $f_\xi [X_s \setminus X_{\alpha(\xi)+1} ] \subseteq \{1\}$ and if $\eta \in [0,\kappa_{\alpha(\xi)+1}]$, then $$f_\xi(\eta,\alpha(\xi) +1) := \begin{cases} 1, & \text{if} \ \eta=\xi, \\
0, & \text{if} \ \eta\neq \xi.
   \end{cases}$$ Also, let $x_\xi := (\xi,\alpha(\xi)+1)$ and $y_\xi := (\kappa_{\alpha(\xi)+1} ,\alpha(\xi)+1)$ whenever $\xi<\kappa$. Consider the family $\{V(f_\xi, \{x_\xi, y_\xi\},1/2) : \xi<\kappa\}$ and suppose, in search of a contradiction, that there are $J \in [\kappa]^{\kappa}$ and $f\in C_p(X_s)$ with $f\in \bigcap\{V(f_\xi, \{x_\xi, y_\xi\}, 1/2) : \xi\in J\}$.
  
Since $f$ is continuous, there is $U_\alpha \subseteq L(\lambda, \kappa_\alpha)$ such that $\kappa_\alpha \in U_\alpha$ and $f[U_\alpha\times \{\alpha\}] = \{f(\kappa_\alpha,\alpha)\}$, provided $\alpha <\cf(\kappa)$. Our goal is to prove that $J\cap J_\alpha \subseteq L(\lambda, \kappa_{\alpha+1})\setminus U_{\alpha+1}$ for any $\alpha <\cf(\kappa)$. Clearly, if $\alpha <\cf(\kappa)$, $J\cap J_\alpha$ is contained in $L(\lambda, \kappa_{\alpha+1})$. Moreover, if we assume the existence of $\xi \in J\cap J_\alpha \cap U_{\alpha+1}$, then the condition $\xi \in U_{\alpha+1}$ implies that $ f(x_\xi) = f(y_\xi)$, while the relation $\xi \in J\cap J_\alpha$ guarantees $$\abs{f(x_\xi) - f_\xi(x_ \xi)} = \abs{f(x_\xi) - 1} < \frac{1}{2} \quad \text{and} \quad \abs{f(y_\xi) - f_\xi(y_ \xi)} = \abs{f(y_\xi)} < \frac{1}{2};$$ consequently, $f(x_\xi) \neq f(y_\xi)$, which is absurd.

Thus, since for each $\alpha<\cf(\kappa)$ it happens that $J\cap J_\alpha \subseteq L(\lambda, \kappa_{\alpha+1})\setminus U_{\alpha+1}$, then $ \abs{J\cap J_\alpha}\leq \abs{L(\lambda, \kappa_{\alpha+1})\setminus U_{\alpha+1}} < \lambda$. Finally, we obtain $\kappa = |J| = |\bigcup_{\alpha<\cf(\kappa)} J\cap J_\alpha| \leq \cf(\kappa) \cdot \sup\{|J\cap J_\alpha| : \alpha<\cf(\kappa)\} \leq \cf(\kappa)\cdot \lambda = \lambda<\kappa$; a contradiction. Hence, $\kappa$ is not a caliber for $C_p(X_s)$.
\end{proof}

\end{document}